\magnification 1200

  \input amssym


  \font \bbfive		= bbm5
  \font \bbseven	= bbm7
  \font \bbten		= bbm10
  \font \eightbf	= cmbx8
  \font \eighti		= cmmi8 \skewchar \eighti = '177
  \font \eightit	= cmti8
  \font \eightrm	= cmr8
  \font \eightsl	= cmsl8
  \font \eightsy	= cmsy8 \skewchar \eightsy = '60
  \font \eighttt	= cmtt8 \hyphenchar\eighttt = -1

  \font \sixi		= cmmi6 \skewchar \sixi = '177
  \font \sixrm		= cmr6
  \font \sixsy		= cmsy6 \skewchar \sixsy = '60
  \font \tensc		= cmcsc10

  \font \titlefont	= cmbx12
  \scriptfont \bffam	= \bbseven
  \scriptscriptfont \bffam = \bbfive
  \textfont \bffam	= \bbten

  \newskip \ttglue

  \def \eightpoint {\def \rm {\fam 0 \eightrm }\relax
  \textfont 0= \eightrm
  \scriptfont 0 = \sixrm \scriptscriptfont 0 = \fiverm
  \textfont 1 = \eighti
  \scriptfont 1 = \sixi \scriptscriptfont 1 = \fivei
  \textfont 2 = \eightsy
  \scriptfont 2 = \sixsy \scriptscriptfont 2 = \fivesy
  \textfont 3 = \tenex
  \scriptfont 3 = \tenex \scriptscriptfont 3 = \tenex
  \def \it {\fam \itfam \eightit }\relax
  \textfont \itfam = \eightit
  \def \sl {\fam \slfam \eightsl }\relax
  \textfont \slfam = \eightsl
  \def \bf {\fam \bffam \eightbf }\relax
  \textfont \bffam = \bbseven
  \scriptfont \bffam = \bbfive
  \scriptscriptfont \bffam = \bbfive
  \def \tt {\fam \ttfam \eighttt }\relax
  \textfont \ttfam = \eighttt
  \tt \ttglue = .5em plus.25em minus.15em
  \normalbaselineskip = 9pt
  \def \MF {{\manual opqr}\-{\manual stuq}}\relax
  \let \sc = \sixrm
  \let \big = \eightbig
  \setbox \strutbox = \hbox {\vrule height7pt depth2pt width0pt}\relax
  \normalbaselines \rm }


  \def \ifundef #1{\expandafter \ifx \csname #1\endcsname \relax }


  \newcount \secno \secno = 0
  \newcount \stno \stno = 0
  \newcount \eqcntr \eqcntr = 0

  \def \track #1#2#3{\ifundef {#1}\else \hbox {\sixrm [#2\string #3] }\fi }

  \def \advseqnumbering {\global \advance \stno by 1 \global \eqcntr =0}

  \def \current {\number \secno \ifnum \number \stno = 0 \else .\number \stno \fi }

  \def \laberr #1{\message {*** RELABEL CHECKED FALSE for #1 ***}
      RELABEL CHECKED FALSE FOR #1, EXITING.
      \end }

  \def \deflabel #1#2{%
    \ifundef {#1}%
      \global \expandafter
      \edef \csname #1\endcsname {#2}%
    \else
      \edef \deflabelaux {\expandafter \csname #1\endcsname }%
      \edef \deflabelbux {#2}%
      \ifx \deflabelaux \deflabelbux \else \laberr {#1=(\deflabelaux )=(\deflabelbux )} \fi
      \fi
    \track {showlabel}{*}{#1}}

  \def \eqmark #1 {\advseqnumbering
    \eqno {(\current )}
    \deflabel {#1}{\current }}

  \def \subeqmark #1 {\global \advance \eqcntr by 1
    \edef \subeqmarkaux {\current .\number \eqcntr }
    \eqno {(\subeqmarkaux )}
    \deflabel {#1}{\subeqmarkaux }}

  \def \label #1 {\deflabel {#1}{\current }}
  \def \lcite #1{(#1\track {showlcit}{$\bullet $}{#1})}


  \catcode `\@=11
  \def \c@itrk #1{{\bf #1}\track {showcitations}{\#}{#1}} 
  \def \c@ite #1{{\rm [\c@itrk{#1}]}}
  \def \sc@ite [#1]#2{{\rm [\c@itrk{#2}\hskip 0.7pt:\hskip 2pt #1]}}
  \def \du@lcite {\if \pe@k [\expandafter \sc@ite \else \expandafter \c@ite \fi }
  \def \cite {\futurelet \pe@k \du@lcite }
  \catcode `\@=12


  \newcount \bibno \bibno = 0
  \newcount \bibtype \bibtype = 0 


  \def \newbib #1#2{\ifcase \number \bibtype
	\global \advance \bibno by 1 \edef #1{\number \bibno }\or
	\edef #1{#2}\or
	\edef #1{\string #1}\fi }

  \def \bibitem #1#2#3#4{\smallbreak \item {[#1]} #2, ``#3'', #4.}

  \def \references {\begingroup \bigbreak \eightpoint
    \centerline {\tensc References}
    \nobreak \medskip \frenchspacing }


  \def \Headlines #1#2{\nopagenumbers
    \headline {\ifnum \pageno = 1 \hfil
    \else \ifodd \pageno \tensc \hfil \lcase {#1} \hfil \folio
    \else \tensc \folio \hfil \lcase {#2} \hfil
    \fi \fi }}

  \def \Title {\centerline {\titlefont \ucase {\titletext }}
    \ifundef {titletextContinued}\else
      \smallskip
      \centerline {\titlefont \ucase {\titletextContinued }}
      \fi }

  \long \def \Quote #1\endQuote {\begingroup \leftskip 35pt \rightskip 35pt
\parindent 17pt \eightpoint #1\par \endgroup }
  \long \def \Abstract #1\endAbstract {\bigskip \Quote \noindent #1\endQuote }
  
  \def \Authors #1{\bigskip \centerline {\tensc #1}}
  \def \Note #1{\footnote {}{\eightpoint #1}}
  \def \Date #1 {\Note {\it Date: #1.}}


  \def \ucase #1{\edef \auxvar {\uppercase {#1}}\auxvar }
  \def \lcase #1{\edef \auxvar {\lowercase {#1}}\auxvar }

  \def \section #1 \par {\global \advance \secno by 1 \stno = 0
    \bigbreak \noindent {\bf \number \secno .\enspace #1.}
    \nobreak \medskip \noindent }

  \def \state #1 #2\par {\medbreak \noindent \advseqnumbering {\bf \current .\enspace #1.\enspace \sl #2\par }\medbreak }
  \def \definition #1\par {\state Definition \rm #1\par }

  \long \def \Proof #1\endProof {\medbreak \noindent {\it Proof.\enspace }#1
\ifmmode \eqno \endproofmarker $$ \else \hfill $\endproofmarker $ \looseness = -1 \fi \medbreak }

  \def \$#1{#1 $$$$ #1}
  \def \explain #1#2{\mathrel {\buildrel \hbox {\sixrm \lcite{#1}} \over #2}}
  \def \=#1{\explain {#1}{=}}

  \def \pilar #1{\vrule height #1 width 0pt}
  \def \stake #1{\vrule depth  #1 width 0pt}

  \newcount \footno \footno = 1
  \newcount \halffootno \footno = 1
  \def \footcntr {\global \advance \footno by 1
  \halffootno =\footno
  \divide \halffootno by 2
  $^{\number \halffootno }$}
  \def \fn #1{\footnote {\footcntr }{\eightpoint #1\par }}




  \def \Item #1{\smallskip \item {{\rm #1}}}
  \newcount \zitemno \zitemno = 0

  \def \izitem {\global \zitemno = 0}
  \def \zitemplus {\global \advance \zitemno by 1 \relax }
  \def \rzitem {\romannumeral \zitemno }
  \def \rzitemplus {\zitemplus \rzitem } 
  \def \zitem {\Item {{\rm (\rzitemplus )}}}
  \def \Zitem {\izitem \zitem }
  \def \zitemmark #1 {\deflabel {#1}{\rzitem }}

  \newcount \nitemno \nitemno = 0
  
  \def \nitem {\global \advance \nitemno by 1 \Item {{\rm (\number \nitemno )}}}

  \newcount \aitemno \aitemno = -1
  \def \boxlet #1{\hbox to 6.5pt{\hfill #1\hfill }}
  \def \iaitem {\aitemno = -1}
  \def \aitemconv {\ifcase \aitemno a\or b\or c\or d\or e\or f\or g\or
h\or i\or j\or k\or l\or m\or n\or o\or p\or q\or r\or s\or t\or u\or
v\or w\or x\or y\or z\else zzz\fi }
  \def \aitem {\global \advance \aitemno by 1\Item {(\boxlet \aitemconv )}}
  \def \aitemmark #1 {\deflabel {#1}{\aitemconv }}


  \font \mf =cmex10
  \def \union {\mathop {\raise 9pt \hbox {\mf S}}\limits }
  \def \inters {\mathop {\raise 9pt \hbox {\mf T}}\limits }

  \def \<{\left \langle \vrule width 0pt depth 0pt height 8pt }
  \def \>{\right \rangle }
  \def \({\big (}
  \def \){\big )}
  \def \ds {\displaystyle }
  \def \and {\hbox {,\quad and \quad }}

  \def \imply {\kern 7pt \Rightarrow \kern 7pt}
  \def \for #1{,\quad \forall \,#1}
  \def \endproofmarker {\square } 
  \def \"#1{{\it #1}\/} 
  \def \inv {^{-1}}
  \def \*{\otimes }
  \def \caldef #1{\global \expandafter \edef \csname #1\endcsname {{\cal #1}}}
  \def \bfdef #1{\global \expandafter \edef \csname #1\endcsname {{\bf #1}}}
  \bfdef N \bfdef Z \bfdef C \bfdef R


  \def \Caixa #1{\setbox 1=\hbox {$#1$\kern 1pt}\global \edef \tamcaixa {\the \wd 1}\box 1}
  \def \caixa #1{\hbox to \tamcaixa {$#1$\hfil }}



  \catcode `\@=11

  \def \overparenOnefill {$\m@th
  \setbox 0=\hbox {$\braceld $}%
  \braceld \leaders \vrule height\ht 0 depth0pt\hfill
  \leaders \vrule height\ht 0 depth0pt\hfill \bracerd $}

  \def \overparenOne #1{\mathop {\vbox {\m@th\ialign {##\crcr \noalign {\kern -1pt}
  \overparenOnefill \crcr \noalign {\kern 3pt\nointerlineskip }
  $\hfil \displaystyle {#1}\hfil $\crcr }}}\limits }

  \def \overparenTwofill {$\m@th
  \lower 0.3pt \hbox {$\braceld $}
  \leaders \vrule depth 0pt height1pt \hfill
  \lower 0.3pt \hbox {$\bracerd $}$}

  \def \overparenTwo #1{\mathop {\vbox {\ialign {##\crcr \noalign {\kern -1pt}
  \overparenTwofill \crcr \noalign {\kern 3pt\nointerlineskip }
  $\hfil \displaystyle {#1}\hfil $\crcr }}}\limits }

  \catcode `\@=12


  %
  \def \adverbOne {r}
  \def \adverbTwo {esidual}
  \def \adjective {free}
  \def \essfree {{\adverbOne }{\adverbTwo }ly \adjective }
  \def \essfreeness {{\adverbOne }{\adverbTwo } \adjective ness}
  \def \Essfreeness {\ucase {\adverbOne }{\adverbTwo } \adjective ness}

  \def \destaca #1{\underline {#1}}   \def \destaca #1{#1}	
  \def \src {d}	
  \def \ran {r}	
  \def \vr {x}	
  \def \vro {y}	
  \def \ed {e}    
  \def \oed {f}	
  \def \eproj {e} 
  \def \s {s}  	
  \def \auto {\destaca {\sigma }}	
  \def \urep {\tilde \pi }  
  \def \ts {\tilde \s }   
  \def \tp {\tilde p}   
  \def \tu {\tilde u}   
  
  \def \Gpd {{\cal G}\tight (\SGE )}
  \def \Gtwo {\Gpd ^{(2)}}
  \def \q {\check }   \def \q {\breve }
  \def \corona {{\q G}}
  \def \lag {\ell }
  \def \O {{\cal O}}

  \def \g {{\bf g}}  
  \def \TEL {\G_{G,E}}
  \def \cyl #1{Z(#1)}

  \def \germ #1#2#3#4{\big [#1,#2,#3;\,#4\big ]}
  \def \trunc #1#2{#1|_{#2}}

  \def \Lin {{\cal L}}

  %
  \def \G {{\cal G}}
  \def \tight {_{\rm tight}}

  \def \Data {G,E,\varphi }
  \def \OGE {{\cal O}_{G,E}}
  
  \def \OAB {{\cal O}_{A,B}}
  \def \S {{\cal S}}
  \def \T {{\cal T}}
  \def \SGE {\S_{G,E}}
  \def \SE {\S_E}
  \def \ISL {{\cal E}}
     
    \def \fix {\smallskip \noindent $\blacktriangleright $\kern 12pt}

  %
  \input pictex

  %
  \def \Subjclass #1#2{\footnote {\null }{\eightrm #1 \eightsl Mathematics Subject Classification:  \eightrm #2.}}
  \def \text #1{\hbox {#1}}
  \def \bool #1{[{\scriptstyle #1}]\,}
  \def \equal #1#2{\bool {#1=#2}}
  \def \? {\vrule width 34pt}

  \def \equationmark #1 {\ifundef {InsideBlock}
	  \advseqnumbering
	  \eqno {(\current )}
	  \deflabel {#1}{\current }
	\else
	  \global \advance \eqcntr by 1
	  \edef \subeqmarkaux {\current .\number \eqcntr }
	  \eqno {(\subeqmarkaux )}
	  \deflabel {#1}{\subeqmarkaux }
	\fi }

  \long \def \Proof #1\endProof {\begingroup \def \InsideBlock {} \medbreak \noindent {\it Proof.\enspace }#1
\ifmmode \eqno \endproofmarker $$ \else \hfill $\endproofmarker $ \looseness = -1 \fi \medbreak \endgroup }

  \def \state #1 #2\par {\begingroup \def \InsideBlock {} \medbreak \noindent \advseqnumbering {\bf \current .\enspace
#1.\enspace \sl #2\par }\medbreak \endgroup }

  %
  \input miniltx \makeatletter \def \Gin @driver{pdftex.def} \input color.sty \resetatcatcode

  %
  \def \ref #1{\sysref #1\endSysRef }
  \def \sysref #1:#2\endSysRef {[[#2]]}

  %

  \newbib \BO {BO}
  \newbib \actions {E}
  \newbib \EP {EP}
  \newbib \ExelVesshik {EV}
  \newbib \Grig {G}
  \newbib \GS {GS}
  \newbib \KatsuraOne {K}
  \newbib \KPRR {KPRR}
  \newbib \Lawson {L}
  \newbib \NekraJO {N}
  \newbib \NC {N2}
  \newbib \pat {P}
  \newbib \Pedersen {Pe}
  \newbib \Pimnsner {Pi}
  \newbib \Raeburn {R}


  \def \titletext {Graphs, groups and self-similarity}

  \Headlines {\titletext } {R.~Exel and E.~Pardo}

  \Title

  \Authors {Ruy Exel and Enrique Pardo}

  \Date {4 July 2013}

  \Subjclass {2010}{46L05, 46L55}

  \Note {\it Key words and phrases: \rm Kirchberg algebra, Katsura algebra, tight representation, inverse semigroup,
groupoid, groupoid C*-algebra.}

  \Note {The first-named author was partially supported by CNPq. The second-named author was partially supported by PAI
III grants FQM-298 and P11-FQM-7156 of the Junta de Andaluc\'{\i}a, by the DGI-MICINN and European Regional Development
Fund, jointly, through Project MTM2011-28992-C02-02 and by 2009 SGR 1389 grant of the Comissionat per Universitats i
Recerca de la Generalitat de Catalunya.}

  \Abstract
  We study a family of C*-algebras generalizing both Katsura algebras and certain algebras introduced by Nekrashevych in
terms of self-similar groups.
  \endAbstract

\section Introduction

The purpose of this paper is to give a unified treatment to two classes of C*-algebras which have been studied in the past
few  years from rather different points of view, namely
  Katsura's algebras \cite{\KatsuraOne},
  and certain algebras constructed by Nekrashevych \cite{\NekraJO}, \cite{\NC}  from  self-similar groups.

The realization that these classes are indeed closely related, as well as the fact that they could be given
a unified treatment, came to our mind as a result of our earlier attempt  \cite{\EP} to  understand Katsura's
algebras $\O_{A,B}$ from the point of view of inverse semigroups. The fact, proven by Katsura in \cite{\KatsuraOne},
that all Kirchberg algebras may be described in terms of  his $\O_{A,B}$ was, in turn, a strong motivation for that endeavor.

While studying $\O_{A,B}$, it slowly became clear to us that the two matricial parameters $A$ and $B$ play very
different roles.  The reader acquainted with Katsura's work will easily recognize that the matrix $A$ is destined to be
viewed as the edge matrix of a graph,
  but it took us much longer to realize that $B$ should be thought of as providing parameters for an action of the group
${\bf Z}$ on the graph given by $A$.  In trying to understand these different roles, some interesting arithmetic popped up
sparking a connection with the work done by Nekrashevych \cite{\NC} on the C*-algebra $\O_{(G,X)}$ associated to a
self-similar group $(G,X)$.

While Nekrashevych's algebras contain a Cuntz algebra, Katsura's algebras contain a graph C*-algebra.   This fact alone
ought to  be considered as a hint that self-similar groups lie in a much bigger class, where the group action takes place  on the path
space of a graph, rather than on a rooted tree (which, incidentally, is the path space of a bouquet of circles).

One of the first important applications of the idea of self-similarity in group theory is in constructing groups with
exotic  properties \cite{\Grig}, \cite{\GS}.  Many of these are defined as subgroups of the group of all automorphisms of a
tree.  Having been born from automorphisms, it is natural that the theory of self-similar groups generally assumes that
the group acts \"{faithfully} on its tree (see, e.g. \cite[Definition 2.1]{\NC}).

However, based on the example provided by Katsura's algebras, we decided that perhaps it is best to  view the group on
its own, the action being an extra ingredient.

The main idea behind self-similar groups, namely the equation
  $$
  g(xw) = yh(w)
  \equationmark SSimilerity
  $$
  appearing in \cite[Definition 2.1]{\NC}, and the subsequent notion   of \"{restriction}, namely
  $$
  g|_x:= h,
  $$
  depend on faithfulness, since otherwise the group element $h$ appearing in \lcite{\SSimilerity} would not be unique
and therefore will not be well
defined as a function of $g$ and $x$.  Working with non-faithful group actions we were forced to postulate a
functional dependence
  $$
  h = \varphi(g,x),
  $$
  and we were surprised to find that the natural properties expected of $\varphi$ are  that of a group cocycle.

To be precise, the ingredients needed in our generalization of self-similar groups are:  a countable discrete group  $G$, an action
  $$
  G\times E\to E
  $$
  of $G$ on a finite graph
  $E = (E^0, E^1, \ran , \src )$, and
  a one-cocycle
  $$
  \varphi:G\times E^1 \to G
  $$
  for  the action of $G$ on the edges of $E$.

Starting with this data (satisfying a few other natural axioms) we construct an action of $G$ on the space of finite paths
$E^*$ which satisfy the ``self-similarity''  equation
  $$
  g(\alpha\beta) = (g\alpha)\big(\varphi(g,\alpha)\beta\big)
  \for g\in E \for \alpha,\beta\in E^*.
  $$

Adopting a philosophy  similar to that  embraced by Katsura and Nekrashevych, we define a C*-algebra, denoted
  $$
  \OGE,
  $$
  in terms of generators and relations inspired by the above group action.  The study of $\OGE$ is, thus, the purpose  of this paper.

Given a  self-similar group $(G,X)$, if we consider $X$ as the set of edges of a graph with a single
vertex, and if we define $\varphi(g,x) = g|_x$, then our $\OGE$  coincides with Nekrashevych's $\O_{(G,X)}$.

On the other hand, if we are given two integer $N\times N$ matrices $A$ and $B$, we may form a graph $E$ with vertex set
$E^0=\{1,2,\ldots,N\}$ and with $A_{i,j}$ edges from vertex $i$ to vertex $j$.  We may then use $B$ to define an action of ${\bf Z}$
on $E$, by fixing all vertexes and acting on the set of edges as follows:
  denote the edges in $E$ from $i$ to $j$ by $e_{i,j,n}$, where $0\leq n<A_{i,j}$.  Given $m\in{\bf Z}$, we perform the
Euclidean division of $mB_{i,j}+n$ by $A_{i,j}$, say
  $$
  mB_{i,j}+n=\hat k A_{i,j} + \hat {n}
  $$
  with  $0\leq\hat {n}<A_{i,j}$, and put
  $$
  \auto_m(e_{i,j,n}) = e_{i,j,\hat n},
  $$
  so that the group element $m$ permutes the $A_{i,j}$ edges from $i$ to $j$ in the same way that addition by
$mB_{i,j}$, modulo $A_{i,j}$, permutes the integers $\{0,1,\ldots,A_{i,j}-1\}$.

The quotient $\hat k$ on the above Euclidean division also plays an important role, being used in the definition of the
cocycle
  $$
  \varphi (m, e_{i,j,n})= \hat k.
  $$

In possession of the graph, the action of ${\bf Z}$, and the cocycle $\varphi$, we apply our construction and find that $\OGE$ is
isomorphic to Katsura's $\OAB$.

So, both Nekrashevych's and Katsura's algebras become special cases of our construction.  We therefore believe
that  the project of
studying such  group actions on  path spaces as well as the corresponding algebras  is of great importance.

Taking the first few steps we have been able to describe $\OGE$ as the C*-algebra of an \'etale groupoid $\TEL$, whose
construction is remarkably similar to the groupoid associated to the relation of ``tail equivalence with
lag'' on the path space,  as described by Kumjian, Pask, Raeburn and Renault in \cite{\KPRR}.

The first similarity  is that our groupoid $\TEL$  has the exact same unit space as the corresponding graph groupoid,
namely the infinite path space.  The second, and most surprising similarity is that $\TEL$ is also described by
a \"{lag} function, except that the values of the lag are  not integer numbers, as in \cite{\KPRR}, but lie in a
slightly more complicated group, the semi-direct product of the corona group of $G$ by the right shift automorphism (see
below for precise definitions).

The techniques we use to give $\OGE$ a groupoid model bear heavily on the theory of tight representations of inverse
semigroups developed by the first named author in \cite{\actions}.  In particular, from our initial data we construct an
abstract inverse semigroup $\SGE$ and show that $\OGE$ is the universal C*-algebra for tight representations of $\SGE$.

As a second step we again take inspiration from Nekrashevych
\cite{\NekraJO} and give a description of $\OGE$ as a Cuntz Pimnsner algebra for a very natural correspondence $M$ over
the algebra
  $$
  C(E^0) \ifundef {rtimes} \times \else \rtimes \fi G.
  $$
  As a result we are able to prove that $\OGE$ is nuclear when the $G$ is amenable.

We would like to stress that, like Nekrashevych \cite[Theorem 5.1]{\NC}, our groupoid $\TEL$ is constructed as a
groupoid of germs.  However, departing from Nekrashevych's techniques, we use Patterson's \cite{\pat} notion of
``germs'', rather than the one employed in \cite[Section 5]{\NC}.  While agreeing in many cases, the former has a much
better chance of producing Hausdorff groupoids and, in our case, we may give a precise characterization of Hausdorffness
in terms of a property we call {\essfreeness} (see below for the precise definition).

Part of this work was done during a visit of the second named author to the Departamento de Matem\'atica da
Universidade Federal de Santa Catarina (Florian\'opolis, Brasil) and he would like to express his thanks to the host center
for its warm hospitality.

\section Groups acting on graphs

  Let $E = (E^0, E^1, \ran , \src )$ be a directed graph, where $E^0$ denotes the set of \"{vertexes}, $E^1$ is the set of
\"{edges}, $\ran $ is the \"{range} map, and $\src $ is the \"{source}, or \"{domain} map.

By definition, a \"{source} in $E$ is a vertex $\vr\in E^0 $, for which  $\ran ^{-1}(\vr )=\ifundef {varnothing} \emptyset \else \varnothing \fi$.
  Thus, when we say that a graph has \"{no sources}, we mean that $\ran ^{-1}(\vr )\neq\ifundef {varnothing} \emptyset
\else \varnothing \fi$, for all $\vr \in E^0$.

  By an \"{automorphism} of $E$ we shall mean a bijective map
  $$
  \auto : E^0 \mathop {\dot \cup} E^1 \to E^0 \mathop {\dot \cup} E^1
  $$
  such that
  $\auto (E^i)\subseteq E^i$, for $i = 0,1$,
  and moreover such that
   $\ran \circ\auto = \auto \circ\ran $, and $\src \circ\auto = \auto \circ\src $, on $E^1$.
  It is evident that the collection of all automorphisms of  $E$ forms a group under composition.

  By an action of a group $G$ on a graph $E$  we shall mean a group homomorphism from $G$ to the group of all
automorphisms of $E$.

If $X$ is any set, and if $\auto $  is an action of a group $G$ on $X$, we shall say that a map
  $$
  \varphi:G\times X \to G
  $$
  is a  \"{one-cocycle} for $\auto $, when
  $$
  \varphi(gh, x) =
  \varphi\big(g,\auto_h(x)\big)\varphi(h,x),
  \equationmark CocycleId
  $$
  for all $g,h \in G$, and all $x \in X$.
Plugging $g = h = 1$ above we see that
  $$
  \varphi(1,x) = 1,
  \equationmark CocycleAtOne
  $$
  for every $x$.

\state {Standing Hypothesis} \label StandingHyp
  \rm Throughout this work we shall let
    $G$ be a countable discrete group,
    $E$ be a finite graph with no sources,
    $\auto $ be an action of $G$ on $E$,
    and
  $$
  \varphi:G\times E^1 \to G
  $$
  be a one-cocycle for the restriction of  $\auto $  to $E^1$, which moreover satisfies
  $$
  \auto_{\varphi(g,\ed )}(\vr ) = \auto_g(\vr )
  \for g \in G \for \ed \in E^1 \for \vr \in E^0.
  \equationmark ActionOfCocycleOnVertex
  $$

The assumptions that $E$ is finite and has no sources will in fact  only be used in the next section and it could probably be removed by
using well known graph C*-algebra techniques.

  By a \"{path} in $E$ of \"{length} $n\geq1$ we shall mean any finite sequence of the form
  $$
  \alpha = \alpha_1\alpha_2\ldots\alpha_n,
  $$
  where $\alpha_i \in E^1$, and $\src (\alpha_i) = \ran (\alpha_{i+1})$, for all $i$ (this is the usual convention when  treating graphs from a
categorical point of view, in which functions compose from right to left).  The \"{range} of $\alpha$ is defined by
  $$
  \ran (\alpha) = \ran (\alpha_1),
  $$
  while the \"{source} of $\alpha$ is defined by
  $$
  \src (\alpha) = \src (\alpha_n).
  $$

  A vertex $\vr \in E^0$ is  considered to be a path of
length zero, in which case we set $\ran (\vr ) = \src (\vr ) = \vr $.

  For every integer $n\geq0$ we denote by $E^n$ the set of all paths in $E$ of length $n$ (this being consistent with the
already introduced notations for $E^0$ and $E^1$).  Finally, we denote by $E^*$ the sets of all finite paths, and by
$E^{\leq n}$ the set of all paths of length at most $n$, namely
  $$
  E^* = \bigcup_{k\geq0}E^k
  \and
  E^{\leq n} = \bigcup_{k=0}^nE^k.
  $$

We will often employ the operation of \"{concatenation} of paths.  That is, if (and only if) $\alpha$ and $\beta$ are paths such that
$\src (\alpha)=\ran (\beta)$, we will denote by $\alpha \beta$ the path obtained by juxtaposing $\alpha$ and $\beta $.

In the special case in which $\alpha $
is a path of length zero, the concatenation $\alpha \beta$ is allowed  if and only if $\alpha =\ran (\beta)$, in which case we set $\alpha \beta =\beta $.
Similarly, when $|\beta |=0$, then $\alpha \beta$ is defined iff $\src (\alpha)=\beta $, and then $\alpha \beta =\alpha $.

We would now like to describe a certain extension of $\auto $ and $\varphi$ to finite paths.

\state Proposition \label extendedaction
  Under the assumptions of  \lcite {\StandingHyp }
  there exists a unique pair $(\auto ^*,\varphi^*)$, formed by an action $\auto ^*$ of $G$ on $E^*$ (viewed simply as a set),
  and a one-cocycle $\varphi^*$ for $\auto ^*$, such that, for every $n\geq0$, every $g \in G$, and every $\vr \in E^0$, one has that:
  \Zitem $\auto ^*_g=\auto_g$, on $E^{\leq1}$,
  \zitem $\varphi^*(g,\vr ) = g$,  \zitemmark CocZero
  \zitem $\varphi^* = \varphi$, on $G\times E^1$, \zitemmark ExtPhi
  \zitem $\auto ^*_g(E^n)\subseteq E^n$, \zitemmark KeepLength
  \zitem $\ran \circ\auto ^*_g=\auto_g\circ\ran $, on $E^n$, \zitemmark MatchRange
  \zitem $\src \circ\auto ^*_g=\auto_g\circ\src $, on $E^n$, \zitemmark MatchSource
  \zitem $\auto_{\varphi^*(g,\alpha)}(\vr )=\auto_g(\vr )$, for all $\alpha \in E^n$, \zitemmark PhiStarOnVert
  \zitem $\auto ^*_1$ is the
  identity\fn {This is evidently already included in the statement that $\auto ^*$ is an action, but we repeat it here to
aid our proof by induction.}
  on $E^n$, \zitemmark PhiOneId
  \zitem $\auto ^*_g(\alpha \beta) = \auto ^*_g(\alpha)\ \auto ^*_{\varphi^*(g,\alpha)}(\beta)$, provided $\alpha$ and
$\beta$ are finite paths with $\alpha \beta \in E^n$, \zitemmark MainConcat
  \zitem $\varphi^*(g, \alpha \beta)=\varphi^*\big(\varphi^*(g,\alpha),\beta \big)$, provided $\alpha$ and $\beta$ are
finite paths with $\alpha \beta \in E^n$. \zitemmark LastItem

\Proof
Initially notice that, once \lcite {\MatchRange }, \lcite {\MatchSource } and \lcite {\PhiStarOnVert } are proved, the
concatenation of the paths
  ``$\auto ^*_g(\alpha)$'' and ``$\auto ^*_{\varphi^*(g,\alpha)}(\beta)$'',
  appearing in \lcite {\MainConcat }, is permitted because
  $$
  \ran \big(\auto ^*_{\varphi^*(g,\alpha)}(\beta)\big) \={\MatchRange }
  \auto_{\varphi^*(g,\alpha)}(\ran (\beta)) \={\PhiStarOnVert }
  \auto_g(\ran (\beta)) =
  \auto_g\big(\src (\alpha)\big) \={\MatchSource }
  \src \big(\auto ^*_g(\alpha)\big).
  $$

For every $g$ in $G$, define $\auto ^*_g$ on $E^{\leq1}$ to coincide with $\auto_g$.  Also, define $\varphi^*$ on $G\times E^{\leq1}$ by
\lcite {\CocZero } and \lcite {\ExtPhi }.  It is then clear that (i--iii) hold and it is easy to see that the remaining
properties (iv--\LastItem ) hold for all $n\leq1$.

We shall complete the definitions of $\auto ^*$ and $\varphi^*$ by induction, so we assume that $m\geq1$,  that
  $$
  \auto ^*_g: E^{\leq m} \to E^{\leq m}
  $$
  is defined  for all $g$ in $G$, that
  $$
  \varphi^*: G\times E^{\leq m} \to G,
  $$
  is defined,
  and that (i--\LastItem ) hold for all $n\leq m$.
We then define
  $$
  \auto ^*_g: E^{m+1} \to E^{m+1}
  $$
  for all $g$ in $G$, and
  $$
  \varphi^*: G\times E^{m+1} \to G,
  $$
  by induction as follows.
  \def \one {\alpha'}\def \two {\alpha''}
  Given $\alpha \in E^{m+1}$, write $\alpha =\one \two $, with $\one \in E^1$, and $\two \in E^m$, and put
  $$
  \auto ^*_g(\alpha) = \auto_g(\one )\auto ^*_{\varphi(g,\one )}(\two )
  \and
  \varphi^*(g, \alpha)=\varphi^*\big(\varphi(g,\one ),\two \big).
  \equationmark InducDef
  $$
  A quick analysis, as done in the first paragraph of this proof, shows that the concatenation of ``$\auto_g(\one )$'' and
``$\auto ^*_{\varphi(g,\one )}(\two )\stake {7pt}$'', appearing above, is permitted.  We next verify  (iv--\LastItem ), substituting $m+1$ for $n$.

We have that the length of $\auto ^*_g(\alpha)$, as defined above, is clearly $1+m$, thus proving (iv).   With respect to
\lcite {\MatchRange } we have that
  $$
  \ran \big(\auto ^*_g(\alpha)\big) =\ran \big(\auto_g(\one )\big) = \auto_g\big(\ran (\one )\big) = \auto_g\big(\ran (\alpha)\big).
  $$
  As for \lcite {\MatchSource }, notice that
  $$
  \src \big(\auto ^*_g(\alpha)\big) = \src \big(\auto ^*_{\varphi(g,\one )}(\two )\big) = \auto_{\varphi(g,\one
)}\big(\src (\two )\big) = \auto_g\big(\src (\two )\big) = \auto_g\big(\src (\alpha)\big).
  $$
  Given $\vr \in E^0$, we have that
  $$
  \auto_{\varphi^*(g,\alpha)}(\vr ) =
  \auto_{\varphi^*(\varphi(g,\one ),\two )}(\vr ) =
  \auto_{\varphi(g,\one )}(\vr ) =
  \auto_g(\vr ),
  $$
  taking care of \lcite {\PhiStarOnVert }.

The verification of \lcite {\PhiOneId } is done as follows:  for $\alpha =\alpha'\alpha''$, as in \lcite {\InducDef }, one has
  $$
  \auto ^*_1(\alpha) =   \auto ^*_1(\alpha'\alpha'') =
  \auto _1(\one )\auto ^*_{\varphi(1,\one )}(\two ) \={\CocycleAtOne }
  \auto _1(\one )\auto ^*_1(\two ) = \one \two =\alpha .
  $$

In order to prove \lcite {\MainConcat }, pick paths $\alpha$ in $E^k$ and $\beta$ in
$E^l$, where $k+l=m+1$, and such that $\src (\alpha)=\ran (\beta)$.

We leave it for the reader to verify \lcite {\MainConcat } in the easy case in which $k=0$, that is, when $\alpha$ is a vertex.
The case $k=1$ is also easy as it is nothing but the  definition of $\auto ^*_g$ given in \lcite {\InducDef }.  So we
may assume that $k\geq2$.

  Writing
$\alpha =\one \two $, with $\one \in E^1$, and $\two \in E^{k-1}$, we then have that $\alpha \beta  = \one \two \beta $, and hence, by definition,
  $$
  \auto ^*_g(\alpha \beta) = \auto_g(\one )\auto ^*_{\varphi(g,\one )}(\two \beta) =
  \auto_g(\one ) \auto ^*_{\varphi(g,\one )}(\two )\ \auto ^*_{\varphi^*(\varphi(g,\one ),\two )}(\beta) \$=
  \auto_g^*(\one \two )\ \auto ^*_{\varphi^*(g,\one \two )}(\beta).
  $$
  We remark that, in last step above, one should use the induction hypothesis in case $k\leq m$, and the definitions of $\auto ^*$ and
$\varphi^*$, when $k=m+1$.

To verify \lcite {\LastItem } we again pick paths $\alpha$ in $E^k$ and $\beta$ in $E^l$, where $k+l=m+1$, and such that
$\src (\alpha)=\ran (\beta)$.  We once more leave the easy case $k=0$ to the reader and observe that the case $k=1$ follows from the
definition of $\varphi^*$.

We may then suppose that $k\geq2$, so we write $\alpha =\one \two $, with $\one \in E^1$, and $\two \in E^{k-1}$.  Then
  $$
  \varphi^*(g,\alpha \beta) =   \varphi^*(g,\alpha'\alpha''\beta) = \varphi^*\big(\varphi(g,\alpha'),\alpha''\beta \big) =
  \varphi^*\Big (\varphi^*\big(\varphi(g,\alpha'),\alpha''\big),\beta \Big ) \$=
  \varphi^*\Big (\varphi^*(g,\alpha'\alpha''\big),\beta \Big ) =   \varphi^*\Big (\varphi^*(g,\alpha \big),\beta \Big ).
  $$

Let us now prove that $\auto ^*$ is in fact an action of $G$ on $E^n$.  We begin by proving that
$\auto ^*_g\auto ^*_h = \auto ^*_{gh}$ on $E^n$, for every $g$ and $h$ in $G$, which we do by induction on $n$.

This follows immediately from the hypothesis for $n\leq1$, so let us assume that $n\geq2$.  Given $\alpha \in E^n$, write $\alpha =\alpha'\alpha''$,
with $\one \in E^1$, and $\two \in E^{n-1}$. Then
  $$
  \auto ^*_g\big(\auto ^*_h(\alpha)\big) =  \auto ^*_g\big(\auto ^*_h(\alpha'\alpha'')\big) =
  \auto ^*_g\big(  \auto_h(\alpha')\auto_{\varphi(h,\alpha')}(\alpha'')  \big) \$=
  \auto_g\big(  \auto_h(\alpha') \big) \auto ^*_{\varphi(g,\auto_h(\alpha'))}\big(\auto_{\varphi(h,\alpha')}(\alpha'')  \big) =
  \auto_{gh}(\alpha')  \auto ^*_{ \varphi(g,\auto_h(\alpha')) \varphi(h,\alpha')}(\alpha'') \$=
  \auto_{gh}(\alpha')  \auto ^*_{ \varphi(gh,\alpha')}(\alpha'') =
  \auto ^*_{gh}(\alpha'\alpha'') =   \auto ^*_{gh}(\alpha).
  $$

  That $\alpha^*_g$ is bijective on each $E^n$ then
  follows\fn {This is why it is useful to include  \lcite {\PhiOneId } as a separate statement, since we may now use it to
prove bijectivity.}
  from \lcite {\PhiOneId }, so $\alpha^*$ is indeed an action of $G$ on $E^n$.

  Finally, let us show that $\varphi^*$ is a cocycle for $\auto ^*$ on $E^n$.  For this fix $g$ and $h$ in $G$ and let $\alpha
 \in E^n$.  Then, with $\alpha=\alpha'\alpha''$, as before,
  $$
  \varphi^*(gh, \alpha)=
  \varphi^*(gh, \alpha'\alpha'')=
  \varphi^*(\varphi(gh, \alpha'),\alpha'')=
  \varphi^*\Big ( \varphi\big(g, \auto_h(\alpha')\big) \varphi(h, \alpha'),\alpha''\Big )\$=
  \varphi^*\Big (\varphi\big(g, \auto_h(\alpha')\big),\auto ^*_{\varphi(h, \alpha')}(\alpha'')\Big )
\varphi^*\big(\varphi(h, \alpha'),\alpha''\big) =: (\star).
  $$
  On the other hand, focusing on the right-hand-side of \lcite {\CocycleId },
notice that
  $$
  \varphi^*(g,\auto ^*_h(\alpha))\varphi^*(h,\alpha) =
  \varphi^*\big(g,\auto ^*_h(\alpha'\alpha'')\big)\varphi^*(h,\alpha'\alpha'') \$=
  \varphi^*\big(g, \auto_h(\alpha')\auto ^*_{\varphi(h,\alpha')}(\alpha'')\big) \varphi^*\big(\varphi(h,\alpha'),\alpha''\big) \$=
  \def \a {\auto_h(\alpha')} \def \b {\auto ^*_{\varphi(h,\alpha')}(\alpha'')}
  \varphi^*\Big (\varphi\big(g, \a \big),\b \Big ) \varphi^*\big(\varphi(h,\alpha'),\alpha''\big),
  $$
  which coincides with $(\star)$ above.  This concludes the proof.
\endProof

\bigskip

The only action of $G$ on $E^*$ to be considered in this paper is $\auto ^*$ so, from now on, we will adopt the shorthand
notation
  $$
  g\alpha = \auto ^*_g(\alpha).
  $$
  Moreover, since $\varphi^*$ extends $\varphi$, we will drop the star decoration and denote $\varphi^*$ simply as $\varphi$.
  The group law, the cocycle condition,  and  properties
  \lcite {\CocZero , \MatchRange , \MatchSource , \PhiStarOnVert , \MainConcat , \LastItem }
  of \lcite {\extendedaction } may then be rewritten as follows:

\state Equations \label Equacoes For every $g$ and $h$ in $G$, for every $\vr \in E^0$, and for every $\alpha$ and $\beta$ in $E^*$ such
that $\src (\alpha)=\ran (\beta)$, one has that
  \iaitem
  \aitem $(gh)\alpha  = g(h\alpha)$,
  \aitem $\varphi(gh, \alpha) =  \varphi\big(g,h\alpha \big)\varphi(h,\alpha),$
  \Item {(\CocZero )} $\varphi(g,\vr ) = g$,
  \Item {(\MatchRange )} $\ran (g\alpha)=g\ran (\alpha)$,
  \Item {(\MatchSource )} $\src (g\alpha)=g\src (\alpha)$,
  \Item {(\PhiStarOnVert )} $\varphi(g,\alpha)\vr =g\vr $,
  \Item {(\MainConcat )} $g(\alpha \beta) = (g\alpha)\ \varphi(g,\alpha)\beta $,
  \Item {(\LastItem )} $\varphi(g, \alpha \beta)=\varphi\big(\varphi(g,\alpha),\beta \big)$.

It might be worth noticing that if $\varphi(g,\alpha)=1$, then \lcite {\Equacoes .\MainConcat } reads
  ``$g(\alpha \beta) = (g\alpha)\beta$'',
  which may be viewed as an associativity property.  However associativity does not hold in general as $\varphi$ is not always
trivial, and hence parentheses must be used.

On the other hand parentheses are unnecessary in expressions of the form $\alpha g\beta$, when $\alpha,\beta \in E^*$, and $g \in G$, since the
only possible interpretation for this expression is the concatenation of $\alpha$ with $g\beta$.

Another useful  property of $\varphi$ is in order.

\state Proposition \label Inverses
  For every $g \in G$, and every $\alpha \in E^*$, one has that
  $$
  \varphi(g^{-1},\alpha) = \varphi(g,g^{-1}\alpha)^{-1}.
  $$

\Proof We have
  $$
  1 =
  \varphi(1,\alpha) =
  \varphi(gg^{-1},\alpha) =
  \varphi(g,g^{-1}\alpha)\varphi(g^{-1},\alpha),
  $$
  from where the conclusion follows.
  \endProof

\section The universal C*-algebra $\OGE $

As in the above section we fix a graph $E$, an action of a group $G$ on $E$, and a one-cocycle $\varphi$ satisfying
\lcite {\StandingHyp }.

  It is our next goal to build  a C*-algebra from this data but first let us recall the following notion from
\cite {\Raeburn }:

\definition \label DefineCKFamily A Cuntz-Krieger $E$-family consists of a set
  $$
  \{p_\vr : \vr \in E^0\}
  $$
  of mutually orthogonal projections and a set
  $$
  \{\s_\ed : \ed \in E^1 \}
  $$
  of partial isometries satisfying
  \izitem
  \zitem $s_\ed ^* s_\ed = p_{\src (\ed )}$, for every $\ed \in E^1$.    \zitemmark CkTwo
  \zitem
  $\ds
  p_\vr =\kern -7pt \sum_{\ed \in \ran ^{-1}(\vr)}s_\ed s_\ed ^*,
  $
  for every $\vr \in E^0$ for which $\ran ^{-1}(\vr )$ is finite and nonempty. \zitemmark CkOne

\definition \label DefineOGE
  We define $\OGE $ to be the universal unital C*-algebra generated by a set
  $$
  \{p_\vr :  \vr \in E^0\}\cup\{\s_\ed :  \ed \in E^1\} \cup\{u_g :  g \in G\},
  $$
  subject to the following relations:
  \iaitem
  \aitem $\{p_\vr :  \vr \in E^0\}\cup\{\s_\ed :  \ed \in E^1\}$ is a Cuntz-Krieger $E$-family, \aitemmark CKRel
  \aitem the map $u:G\rightarrow \OGE $, defined by the rule $g\mapsto u_g$, is a unitary representation of $G$,
  \aitem $u_g\s_\ed =\s_{g\ed }u_{\varphi(g,\ed )}$, for every $g \in G$, and $\ed \in E^1$,
  \aitem $u_gp_\vr =p_{g\vr }u_g$, for every $g \in G$, and $\vr \in E^0$.

\bigskip Observe that, under our standing assumptions \lcite {\StandingHyp }, for every $\vr \in E^0$ we have that $\ran ^{-1}(\vr )$
is finite and nonempty.  So \lcite {\DefineCKFamily .ii} and
\lcite {\DefineOGE .\CKRel } imply that
  $$
  u_gp_\vr u_g^* =
  \sum_{\ran (\ed )=\vr } u_g\s_\ed \s_\ed ^*u_g^* =
  \sum_{\ran (\ed )=\vr } \s_{g\ed }u_{\varphi(g,\ed )}u_{\varphi(g,\ed )}^*\s_{g\ed }^* \$=
  \sum_{\ran (\ed )=\vr } \s_{g\ed }\s_{g\ed }^* =
  \sum_{\ran (\oed )=g\vr } \s_{\oed }\s_{\oed }^* =
  p_{g\vr },
  $$
  which says that  \lcite {\DefineOGE .d} follows  from the other conditions.  We have nevertheless included it in
\lcite {\DefineOGE} in the belief that our theory may be generalized to graphs with sources.


  %

Our construction generalizes some well known constructions in the literature as we would now like to mention.

\state Example \rm \label Nekrashevych
  Let $(G,X)$ be a self similar group   as in \cite [Definition 2.1]{\NC }.
  We may then consider a graph $E$ having only
one vertex and such that $E^1=X$.  If we define
  $$
  \varphi(g,\vr )=g|_\vr,
  $$
  where, in the terminology of \cite {\NC }, $g|_\vr $ is the  restriction (or section) of $g$ at $\vr $,
  then the triple $(\Data )$ satisfies \lcite {\StandingHyp } and one may show that $\OGE $ is isomorphic to the algebra
${\cal O}_{(G,X)}$ introduced by Nekrashevych in \cite {\NC }.

\state Example \rm \label Katsura As in \cite {\KatsuraOne }, given a positive integer $N$, let $A \in M_N(\Z ^+)$ without
zero rows and let $B \in M_N(\Z )$ be such that
  $$
  A_{i,j}=0 \imply B_{i,j}=0.
  $$
  Consider the graph $E$ with
  $$
  E^0 = \{1,2,\ldots,N\},
  $$
  and whose  adjacency matrix is $A$.  For each pair of vertexes $i,j \in E^0$,
such that $A_{i,j}\ne 0$, denote the set of edges with range $i$ and source $j$  by
  $$
  \{e_{i,j,n} : 0\leq n<A_{i,j}\}.
  $$

Define an action $\auto $ of ${\bf Z}$ on $E$, which is the identity on $E^0$, and which acts on edges as follows: given $m\in{\bf Z}$, and
  $e_{i,j,n} \in E^1$, let $(\hat k,\hat n)$ be the unique pair of integers such that
  $$
  mB_{i,j}+n=\hat k A_{i,j} + \hat {n} \and 0\leq\hat {n}<A_{i,j}.
  $$
  Thus, $\hat k$ is the quotient and $\hat {n}$ is the remainder of the division of
$mB_{i,j}+n$ by $A_{i,j}$.  We then put
  $$
  \auto_m(e_{i,j,n})=e_{i,j,\hat {n}}.
  $$
  In other words,  $\auto_m$ corresponds to the addition of  $mB_{i,j}$ to the variable ``$n$'' of ``$e_{i,j,n}$'', taken modulo $A_{i,j}$.
  In turn, the one-cocycle is defined by
  $$
  \varphi (m, e_{i,j,n})= \hat k.
  $$

  It may be proved without much difficulty that $\OGE $ is isomorphic to Katsura's \cite {\KatsuraOne } algebra $\OAB $,
under an isomorphism sending each $u_m$ to the $m ^{th}$ power of the unitary
  $$
  u:=\sum_{i=1}^Nu_i
  $$
  in $\OAB $,  and sending $\s_{e_{i,j,n}}$ to $\s_{i,j,n}$.

\bigskip
When $N=1$, the relevant graph for Katsura's  algebras is the same as the one we used above in the description of
Nekrashevych's example.  However the former is not a special case of the latter because, contrary to what is required
in \cite {\NC }, the group action might not be faithful.

  %

We now return to the general case of a triple $(\Data )$ satisfying \lcite {\StandingHyp }.  We initially recall the
usual extension of the notation ``$\s_\ed $'' to allow for paths of arbitrary length.


\state Definition \label ExtendToPath
  Given a finite path $\alpha$ in $E^*$, we shall let $\s_\alpha$ denote the element of $\OGE $ given by:
  \izitem
  \zitem when $\alpha=\vr \in E^0$, we let $\s_\alpha = p_\vr $,
  \zitem when $\alpha \in E^1$, then $\s_\alpha$ is already defined above,
  \zitem when $\alpha \in E^n$, with $n>1$, write $\alpha = \alpha'\alpha''$, with $\alpha' \in E^1$, and $\alpha'' \in E^{n-1}$, and set
$\s_\alpha = \s_{\alpha'}\s_{\alpha''}$, by recurrence.

Commutation  relation \lcite {\DefineOGE .c} may then be generalized to finite paths as follows:

\state Lemma \label RuleforFinitePaths
Given $\alpha \in E^*$ and  $g \in G$, one has that
  $$
  u_g\s_\alpha =\s_{g\alpha}u_{\varphi (g,\alpha)}.
  $$

\Proof
  Let $n$ be the length of $\alpha$.  When $n=0,1$, this follows from \lcite {\DefineOGE .d\&c}, respectively.
When $n>1$, write $\alpha = \alpha'\alpha''$, with $\alpha' \in E^1$, and $\alpha'' \in E^{n-1}$.  Using induction, we then have
  $$
  u_g\s_{\alpha} =
  u_g\s_{\alpha'}\s_{\alpha''} =
  \s_{g\alpha'}u_{\varphi(g,\alpha')}\s_{\alpha''} =
  \s_{g\alpha'}\s_{\varphi(g,\alpha')\alpha''} u_{\varphi(\varphi(g,\alpha'),\alpha'')} \$=
  \s_{(g\alpha')\varphi(g,\alpha')\alpha''} u_{\varphi(g,\alpha'\alpha'')} =
  \s_{g(\alpha'\alpha'')} u_{\varphi(g,\alpha'\alpha'')} =
  \s_{g\alpha} u_{\varphi(g,\alpha)}.
  \endProof

Our next result provides a spanning set for $\OGE $.

\state Proposition \label InvSemPicture
  Let
  $$
  \S = \big \{\s_\alpha u_g\s_\beta^* :  \alpha,\beta \in E^*,\ g \in G,\ \src (\alpha)=g  \src (\beta)\big \} \cup\{0\}.
  $$
  Then $\S $ is closed under multiplication and adjoints and its closed linear span coincides with  $\OGE $.

\Proof That $\S $ is closed under adjoints is clear.  With respect to closure under multiplication, let $\s_\alpha
u_g\s_\beta^*$ and $\s_\gamma u_h\s_\delta^*$ be elements of $\S $.

From \lcite {\DefineOGE .\CKRel } we know that $\s_\beta^*\s_\gamma= 0$, unless either $\gamma=\beta\varepsilon$, or $\beta=\gamma\varepsilon$, for some
$\varepsilon \in E^*$.  If $\gamma=\beta\varepsilon$, then
  $$
  \s_\beta^*\s_\gamma = \s_\beta^*\s_{\beta\varepsilon} = \s_\beta^*\s_\beta\s_\varepsilon = \s_\varepsilon,
  $$
  and hence
  $$
  (\s_\alpha u_g\s_\beta^*)(\s_\gamma u_h\s_\delta^*)=
  \s_\alpha u_g\s_\varepsilon u_h\s_\delta^*=
  \s_\alpha\s_{g  \varepsilon}u_{\varphi (g,\varepsilon)}u_h\s_\delta^* =
  \s_{\alpha g  \varepsilon}u_{\varphi (g,\varepsilon)h}\s_\delta^*.
  \equationmark ProdSpan
  $$
  Moreover, since
  $$
  \src (\alpha g  \varepsilon) =
  \src (g\varepsilon) =
  g\src (\varepsilon) =
  \varphi(g,\varepsilon)\src (\varepsilon) =
  \varphi(g,\varepsilon)\src (\gamma) =
  \varphi(g,\varepsilon)h\src (\delta),
  $$
  we deduce that the element appearing in the right-hand-side of \lcite {\ProdSpan } indeed belongs to $\S $.

In the second case, namely if  $\beta=\gamma\varepsilon$, then the adjoint of the term appearing in the left-hand-side of
\lcite {\ProdSpan } is
  $$
  (\s_\delta u_{h\inv }\s_\gamma^*)
  (\s_\beta u_{g\inv }\s_\alpha^*),
  $$
  and the case already dealt with implies that this belongs to $\S $.  The result then follows from the fact that $\S $ is
self-adjoint.

In order to  prove that $\OGE $ coincides with the closed linear span of $\S $, let $A$ denote the latter.  Given that
$\S $ is self-adjoint and closed under multiplication, we see that $A$ is a closed *-subalgebra of $\OGE $.  Since $A$
evidently contains  $\s_\alpha$ for every $\alpha$ in $E^{\leq1}$, and that it also contains $u_g$ for every $g$ in $G$, we
deduce that $A=\OGE $.
\endProof

\section The inverse semigroup $\SGE $

As before, we keep \lcite {\StandingHyp } in force.

In this section we will give an abstract description of the set $\S $ appearing in \lcite {\InvSemPicture } as well as its
multiplication and  adjoint operation.  The goal is to construct an inverse semigroup from which we will later recover $\OGE $.

\definition \label inverseSemigroup
  Over the set
  $$
  \SGE =\big \{ (\alpha,g,\beta) \in E^*\times G\times E^*: \src (\alpha)=g\src (\beta)\big \}\cup\{0\},
  $$
  consider a  binary \"{multiplication} operation defined by
  $$
  (\alpha,g,\beta)  (\gamma,h,\delta) = \left \{\matrix {
  (\alpha g  \varepsilon,\hfil \varphi(g,\varepsilon) h,\hfil \delta), &  \text {if } \gamma=\beta \varepsilon,    \cr \cr
  (\alpha,\  g\varphi(h^{-1},\varepsilon)^{-1},\  \delta h^{-1}  \varepsilon), &   \text {if } \beta=\gamma \varepsilon,  \cr \cr
  0,\hfil \hfil \hfil & \text {otherwise,}
  }\right .
  $$
  and a unary \"{adjoint} operation defined by
  $$
  (\alpha,g,\beta)^*:= (\beta,g^{-1}, \alpha).
  $$
  The subset of $\SGE $ formed by all elements $(\alpha,g,\beta)$, with $g=1$, will be denoted by $\SE $.

It is easy to see that $\SE $ is closed under the above operations, and that it is  isomorphic to the
inverse semigroup generated by the canonical partial isometries in the graph C*-algebra of $E$.

Let us begin with a simple, but useful result:

\state Lemma \label EasyMult
  Given $(\alpha,g,\beta)$ and $(\gamma,h,\delta)$ in $\SGE $, one has
  $$
  \beta=\gamma \imply (\alpha,g,\beta) (\gamma,h,\delta) =  (\alpha,gh,\delta).
  $$

\Proof
  Focusing on the first clause of \lcite {\inverseSemigroup }, write $\gamma=\beta\varepsilon$, with $\varepsilon=\src (\beta)$.  Then
  $$
  (\alpha,g,\beta)  (\gamma,h,\delta)=
  (\alpha g \src (\beta),\  \varphi\big(g,\src (\beta)\big) h,\  \delta) =
  (\alpha  \src (\alpha),\ gh,\  \delta) =   (\alpha,gh,  \delta).
  \endProof

\state Proposition $\SGE $ is an inverse semigroup with zero.

\Proof We leave it for the reader to prove that the above  operations are  well defined and associative.
In order to prove the statement it then suffices \cite [Theorem 1.1.3]{\Lawson} to show that, for all $y,z \in \SGE $, one has that
  \Zitem $yy^*y=y$, and
  \zitem $yy^*$ commutes with $zz^*$.

\bigskip

Given $y=(\alpha,g,\beta) \in \SGE $, we have by the above Lemma that
  $$
  yy^*y =
  (\alpha,g,\beta)(\beta, g^{-1}, \alpha)(\alpha,g,\beta)=
  (\alpha, 1, \alpha)(\alpha,g,\beta) =
  (\alpha,g,\beta) = y,
  $$
  proving (i).
  Notice also that
  $$
  yy^* = (\alpha, 1, \alpha)
  \equationmark IdempotForm
  $$
  is an element of the idempotent semi-lattice of $\SE $, which is a commutative set because $\SE $ is an inverse
semigroup.  Point (ii) above then follows immediately, concluding the proof.
  \endProof

As seen in  \lcite {\IdempotForm }, the idempotent semi-lattice of $\SGE $, henceforth denoted by $\ISL $,
is given by
  $$
  \ISL = \big \{(\alpha,1,\alpha): \alpha \in E^*\big \} \cup \{0\}.
  \equationmark ISLDef
  $$
  Evidently $\ISL $ is also  the idempotent semi-lattice of $\SE $.

  For simplicity, from now on we will adopt the
short-hand  notation
  $$
  \eproj_\alpha = (\alpha,1,\alpha) \for \alpha \in E^*.
  \equationmark DefineEAlpha
  $$

The following is a standard fact in the theory of graph C*-algebras:

\state Proposition \label ProdIdemp If  $\alpha,\beta \in E^*$, then
  $$
  \eproj_\alpha\eproj_\beta = \left \{\matrix {
  \eproj_\alpha, &  \text {if there exists $\gamma$ such that } \alpha=\gamma \beta,    \cr
  \eproj_\beta, &  \text {if there exists $\gamma$ such that } \alpha\gamma=\beta, \pilar {12pt}    \cr
  \hfill 0, & \text {otherwise.}\hfill \pilar {12pt}
  }\right .
  $$

Recall that if $\alpha$ and $\beta$ are in $E^*$, we say that $\alpha \preceq \beta$, if $\alpha$ is a \"{prefix} of $\beta$, i.e.~if there
exists $\gamma \in E^*$, such that $\alpha\gamma = \beta$.  It therefore follows from \lcite {\ProdIdemp } that
  $$
  \eproj_\alpha\leq\eproj_\beta \iff \beta\preceq \alpha.
  \equationmark WrodVsIdemp
  $$

Another easy consequence of \lcite {\ProdIdemp } is that, for any two elements $e,f \in \ISL $, one has that either $e\perp f$,
or $e$ and $f$ are comparable.   In other words
  $$
  e\Cap f \imply e\leq f \text {, \ or \ }  f\leq e.
  \equationmark Compara
  $$
  Recall that, according to \cite [Definition 11.1]{\actions }, we say that  $e$ \"{intersects} $f$, in symbols $e\Cap f$, when $ef\neq0$.

\section {\Essfreeness } and E*-unitarity

Again working under \lcite{\StandingHyp}, suppose we are given $g$ in $G$ and $\alpha$ in $E^*$ such that $g\alpha =
\alpha$, and $\varphi(g,\alpha)=1$.  Then, by \lcite {\Equacoes .\MainConcat }, we have that
  $$
  g(\alpha \beta) = \alpha\beta,
  $$
  whenever $\src (\alpha)=\ran (\beta)$.  Such an element  $g$ therefore acts trivially on a large set of finite words.

Occasionally this will be an annoying feature which we would rather avoid, so we make the following:

\definition \label EssFree
  We will say that $(\Data )$ is \"{\essfree } if, whenever $(g,\ed ) \in G\times E^1$, is  such that $g\ed = \ed $, and $\varphi(g,\ed )=1$, then
$g=1$.

This property may be generalized to finite paths:

\state Proposition \label EssFreePath
  Suppose  that $(\Data )$ is {\essfree } and that $(g,\alpha) \in G\times E^*$, is  such that $g\alpha = \alpha$, and
$\varphi(g,\alpha)=1$, then $g=1$.

\Proof Assume that there is a counter-example $(g,\alpha)$ to the statement, which we assume is minimal
in the sense that $|\alpha|$ is as small as possible.

To be sure, to say that $(g,\alpha)$ is a counter-example is to say that $g\alpha = \alpha$, \
$\varphi(g,\alpha)=1$, and yet  $g\neq1$.

By \lcite {\Equacoes .\CocZero }, $\alpha$ can't be a vertex, and neither can it be an edge, by hypothesis.  So $|\alpha|\geq2$, and we
may then write $\alpha=\beta\gamma$, with $\beta,\gamma\in E^*$, and $|\beta|,|\gamma|<|\alpha|$.  Then
  $$
  \beta\gamma =  \alpha = g\alpha = g(\beta\gamma)  = (g\beta)\varphi(g,\beta)\gamma,
  $$
  whence $\beta=g\beta$, and $\gamma = \varphi(g,\beta)\gamma$, by length considerations.  Should $\varphi(g,\beta)=1$,
the pair $(g,\beta)$ would be a smaller counter-example to the statement, violating the minimality of $\alpha$.  So we
have that $\varphi(g,\beta)\neq1$.  In addition,
  $$
  \varphi\big(\varphi(g,\beta),\gamma\big) = \varphi(g,\beta\gamma) = \varphi(g,\alpha) = 1.
  $$
  It follows that $\big(\varphi(g,\beta),\gamma\big)$ is a counter-example to the statement, violating the minimality of $\alpha$.  This is a
contradiction and hence no counter-example exists whatsoever, concluding the proof.
  \endProof

An apparently stronger version of {\essfreeness } is in order.

\state Proposition \label VarStar
  Suppose that $(\Data )$ is {\essfree }.  Then, for all $g_1,g_2 \in G$, and $\alpha \in E^*$, one has that
  $$
  g_1\alpha=g_2\alpha \quad \wedge \quad \varphi(g_1,\alpha)=\varphi(g_2,\alpha) \imply g_1=g_2.
  $$

\Proof
  Defining $g=g_2^{-1} g_1$, observe that
  $
  g\alpha=\alpha,
  $
  and we claim that $\varphi(g,\alpha) = 1$.
  In fact,
  $$
  \varphi(g,\alpha) =
  \varphi\big(g_2^{-1}g_1,\alpha\big) =
  \varphi\big(g_2^{-1},g_1\alpha\big) \varphi\big(g_1,\alpha\big) \={\Inverses }
  \varphi\big(g_2,g_2^{-1}g_1\alpha\big)^{-1} \varphi(g_1,\alpha) \$=
  \varphi(g_2,\alpha)^{-1} \varphi(g_1,\alpha) =1,
  $$
  so it follows that $g=1$, which is to say that $g_1=g_2$.
  \endProof

Recall that an inverse semigroup $\S $ with zero is called E*-unitary, or sometimes 0-E-unitary \cite[Chapter
9]{\Lawson} if, whenever an element
$s\in\S $ dominates a nonzero idempotent $e$, meaning that
  $se=e$,
  then  $s$ is necessarily also idempotent.

{\Essfreeness } is very closely related to E*-unitary inverse semigroups, as we would now like to show.

\state Proposition \label EssFreeEUnitary
  $(\Data )$ is {\essfree } if and only if $\SGE $ is an E*-unitary inverse semigroup.

\Proof Assuming that $(\Data )$ is {\essfree }, let $s=(\alpha,g,\beta)$ be in $\SGE $, and let $\eproj_\gamma=(\gamma,1,\gamma)$ be a nonzero
idempotent in $\ISL $, such that $\eproj_\gamma\leq s$.  It follows that also
  $$
  \eproj_\gamma\leq s^*s = (\beta,1,\beta),
  $$
  so $\gamma=\beta\varepsilon$, for some $\varepsilon\in E^*$, by  \lcite {\ProdIdemp }.  The relation ``$\eproj_\gamma=s\eproj_\gamma$'' translates into
  $$
  (\gamma,1,\gamma)=(\alpha,g,\beta) (\gamma,1,\gamma) =  \big(\alpha g\varepsilon, \varphi (g,\varepsilon), \gamma\big),
  \equationmark EUniEq
  $$
  so
  $$
  \beta\varepsilon = \gamma = \alpha g\varepsilon,
  $$
  which implies that $\varepsilon=g\varepsilon$, and $\beta=\alpha$.  If we further notice that \lcite {\EUniEq } gives
$\varphi(g,\varepsilon)=1$, we may conclude from \lcite {\EssFreePath } that $g=1$, and hence that
  $$
  s = (\alpha,g,\beta) = (\beta,1,\beta)
  $$
  is an idempotent element.  This shows that $\SGE $ is E*-unitary.

In order to prove the converse, let $(g,\ed )\in G\times E^1$ be such that $g\ed =\ed $, and $\varphi(g,\ed )=1$.  Then
  $$
  (1,g,1)  (\ed ,1,\ed )  = (g\ed ,\varphi(g,\ed ),\ed ) = (\ed ,1,\ed ).
  $$
  Therefore the nonzero element $(1,g,1)$ dominates the idempotent $(\ed ,1,\ed )$ and, assuming that $\SGE $ is E*-unitary,
we conclude that $(1,g,1)$ is idempotent, which is to say that $g=1$.  This proves that $(\Data )$ is {\essfree }.
\endProof

\section Tight representations of $\SGE $

As before, we keep \lcite {\StandingHyp } in force.

It is the main goal of this section to show that $\OGE $ is the
universal C*-algebra for tight representations of $\SGE $.

Recall from \lcite {\ProdIdemp } that $\eproj_\alpha\leq\eproj_{\src (\alpha)}$, for every $\alpha \in E^*$, so we see that  the set
  $$
  \{\eproj_\vr : \vr \in E^0\}
  \equationmark FiniteCover
  $$
  is a \"{cover} \cite [Definition 11.5]{\actions } for $\ISL $.

\state Proposition \label TightRep
  The map
  $$
  \pi : \SGE \to \OGE ,
  $$
  defined by $\pi(0)=0$, and
  $$
  \pi(\alpha,g,\beta) = \s_\alpha u_g\s_\beta^*,
  $$
  is a tight \cite [Definition 13.1]{\actions } representation.


\Proof We leave it for the reader to show that $\pi$ is in fact multiplicative and that it preserves adjoints.

In order to prove that $\pi$ is tight, we shall use the characterization given in \cite [Proposition 11.8]{\actions },
observing that $\pi$ satisfies condition (i) of \cite [Proposition 11.7]{\actions } because, with respect to the cover
\lcite {\FiniteCover }, we have that
  $$
  \bigvee_{\vr \in E^0} \pi(\eproj_\vr ) =   \sum_{\vr \in E^0} \pi(\eproj_\vr ) =   \sum_{\vr \in E^0} p_\vr = 1,
  $$
  by   \lcite {\DefineOGE .\CKRel }.
  So we assume that $\{\eproj_{\alpha_1},\ldots,\eproj_{\alpha_n}\}$ is a cover for a given $\eproj_\beta$, where
$\alpha_1,\ldots\alpha_n,\beta \in E^*$, and we need to show that
  $$
  \bigvee_{i=1}^n\pi(\eproj_{\alpha_i}) \geq \pi(\eproj_\beta).
  \equationmark TightGoal
  $$

In particular, for each $i$, we have that $\eproj_{\alpha_i}\leq\eproj_\beta$, which says that there exists $\gamma_i
\in E^*$ such that $\alpha_i=\beta\gamma_i$.

We shall prove \lcite {\TightGoal } by induction on the variable
  $$
  L = \min_{1\leq i\leq n}|\gamma_i|.
  $$

  If $L=0$, we may pick $i$ such that $|\gamma_i|=0$, and then necessarily $\gamma_i= \src (\beta)$, in which case $\alpha_i=\beta$,
and \lcite {\TightGoal } is trivially true.

Assuming that $L\geq1$, one sees that $\vr := \src (\beta)$ is not a source, meaning that $\ran ^{-1}(\vr )$ is nonempty.  Write
  $$
  \ran ^{-1}(\vr ) = \{\ed _1,\ldots,\ed_k\},
  $$
  and observe that
  $$
  \pi(\eproj_\beta) = \s_\beta\s_\beta^* = \s_\beta p_\vr \s_\beta^* \={\DefineOGE .\CKRel }
  \sum_{j=1}^k\s_\beta \s_{\ed_j}\s_{\ed_j}^*\s_\beta^* =
  \sum_{j=1}^k \pi(\eproj_{\beta \ed_j}).
  $$
  In order to prove \lcite {\TightGoal } it is therefore enough to show that
  $$
  \bigvee_{i=1}^n\pi(\eproj_{\alpha_i}) \geq \pi(\eproj_{\beta \ed_j}),
  \equationmark SmallGoal
  $$
  for all $j=1,\ldots,k$.

Fixing $j$ we claim that $\eproj_{\beta \ed_j}$ is covered by the set
  $$
  Z = \big \{\eproj_{\alpha_i}: 1\leq i\leq n,\ \eproj_{\alpha_i}\leq\eproj_{\beta \ed_j}\big \}.
  $$
  In order to see this let $x \in \ISL $ be a nonzero element such that $x\leq\eproj_{\beta \ed_j}$.   Then $x\leq\eproj_{\beta}$, and so
$x\Cap \eproj_{\alpha_i}$ for  some $i$.
Thus,  to prove the claim it is enough to check that $\eproj_{\alpha_i}$ lies in $Z$.
  Observe that
  $$
  x \eproj_{\beta \ed_j} \eproj_{\alpha_i} =
  x  \eproj_{\alpha_i} \neq 0,
  $$
  which implies that $\eproj_{\beta \ed_j} \Cap \eproj_{\alpha_i}$.

By \lcite {\Compara } we have that
  $\eproj_{\beta \ed_j}$ and $\eproj_{\alpha_i}$ are comparable, so  either
$\beta \ed_j\preceq \alpha_i$ or  $\alpha_i\preceq \beta \ed_j$, by \lcite {\WrodVsIdemp }.
  Since we are under the hypothesis that $L\geq1$, and hence that
  $$
  |\alpha_i| =
  |\beta_i| + |\gamma_i|  \geq
  |\beta|+1 =
  |\beta \ed_j|,
  $$
  we must have that $\beta \ed_j\preceq \alpha_i$, from where we deduce that $\eproj_{\alpha_i}\leq \eproj_{\beta \ed_j}$, proving
our claim.

Employing the induction hypothesis we then deduce that
  $$
  \bigvee_{z \in Z}\pi(z) \geq \pi(\eproj_{\beta \ed_j}),
  $$
  verifying  \lcite {\SmallGoal }, and thus concluding the proof.
  \endProof


We would now  like to prove that the representation $\pi$ above is in fact the \"{universal} tight representation of
$\SGE $.

  \def \urep {\rho}  
  \def \ts {\tilde \s }   
  \def \tp {\tilde p}   
  \def \tu {\tilde u}   

\state Theorem Let $A$ be a unital C*-algebra and let $\urep :\SGE \to A$ be a tight representation.  Then there exists
a unique unital *-homomorphism $\psi:\OGE \to A$, such that the diagram
  \hfill \break
  \vbox {
  \beginpicture
  \setcoordinatesystem units <0.0040truecm, -0.0040truecm> point at 3000 0
  \setplotarea x from -1500 to 1000, y from 500 to 1000
  \put {
    $\matrix {
      \SGE & \buildrel \ds \pi \over \longrightarrow& \OGE \cr \cr
      && \ \Big \downarrow \psi  \cr \cr
      && A \
      }$
    } at 500 500
  \arrow <0.11cm> [0.3,1.2] from 300 400 to 650 650 \put {$\urep $}  at 400 600
  \endpicture
  }
  \break commutes.

\Proof We will initially prove that the elements
  $$
  \matrix {
  \tp_\vr := \urep (\vr ,1,\vr ), \hfill & \forall\vr \in E^0, \hfill \cr \cr
  \ts_\ed := \urep \big(\ed , 1, \src (\ed )\big), \hfill & \forall\ed \in E^1, \hfill \cr \cr
  \tu_g := \ds \sum_{\vr \in E^0}\urep (\vr , g, g^{-1} \vr ), \hfill & \forall g \in G, \hfill }
  $$
  satisfy relations \lcite {\DefineOGE .a--d}.
  Since the $\eproj_\vr $ (defined in \lcite {\DefineEAlpha }) are mutually orthogonal idempotents in $\SGE $, it is clear that
the $\tp_\vr $ are mutually orthogonal projections.  Evidently the $\ts_\ed $ are partial isometries so, in order to check
\lcite {\DefineOGE .a}, we must only verify \lcite {\DefineCKFamily .\CkTwo } and \lcite {\DefineCKFamily .\CkOne }.  With
respect to the former, let $\ed \in E^1$.  Then
  $$
  \ts_\ed ^* \ts_\ed =
  \urep \big(( \src (\ed ), 1, \ed )(\ed , 1, \src (\ed )\big) =
  \urep \big(\src (\ed ), 1, \src (\ed )\big) =
  \tp_{\src (\ed )},
  $$
  proving \lcite {\DefineCKFamily .\CkTwo }.  In order to prove \lcite {\DefineCKFamily .\CkOne }, let $\vr $ be a vertex such
that $\ran ^{-1}(\vr )$ is nonempty and write
  $$
  \ran ^{-1}(\vr ) = \big \{\ed _1,\ldots,\ed_n\big \}.
  $$

  Putting $q_i = (\ed_i,1,\ed_i)$,
  we then claim that the set
  $$
  \big \{q_1,\ldots,q_n\big \}
  $$
  is a cover for $q:= (\vr ,1,\vr )$.  In order to prove this we must show that, if the nonzero idempotent $f$ is dominated by
$q$, then $f\Cap q_i$ for some $i$.

Let $f=(\alpha,1,\alpha)$ by \lcite {\ISLDef } and notice that
  $$
  0 \neq f = fq = (\alpha,1,\alpha)(\vr ,1,\vr ).
  $$
  So $\alpha$ and $\vr $ are comparable, and this can only happen when $\vr =\ran (\alpha)$.  If $|\alpha| = 0$ then necessarily
$\alpha=\vr $, so $f=q$, and it is clear that $f\Cap q_i$ for all $i$.  On the other hand, if $|\alpha|\geq1$, we write
  $$
  \alpha=\alpha'\alpha'',
  $$
  with $\alpha' \in E^1$, so that $\ran (\alpha') = \ran (\alpha)=\vr $, and hence $\alpha'=\ed_i$, for some $i$.
  Therefore
  $$
  fq_i =
  (\alpha,1,\alpha)(\ed_i,1,\ed_i) =
  (\alpha,1,\alpha)(\alpha',1,\alpha') =
  (\alpha,1,\alpha)\neq0,
  $$
  so $f\Cap q_i$, proving the claim.  Since $\urep $ is a tight representation, we deduce that
  $$
  \urep (q) = \bigvee_{i=1}^n\urep (q_i),
  $$
  but since the $q_i$ are easily seen to be pairwise orthogonal, their supremum coincides with their sum, whence
  $$
  \tp_\vr = \urep (q) = \sum_{i=1}^n\urep (q_i) =  \sum_{i=1}^n\urep (\ed_i,1,\ed_i) \$=
  \sum_{i=1}^n \urep \big((\ed_i,1,\src (\ed_i))\ (\src (\ed_i),1,\ed_i) \big) =
  \sum_{i=1}^n \ts_{\ed_i}\ts_{\ed_i}^*,
  $$
  thus verifying \lcite {\DefineCKFamily .\CkOne }, and hence proving \lcite {\DefineOGE .a}.

With respect to \lcite {\DefineOGE .b}, let us first prove that $\tu _1 = 1$.  Considering the subsets of $\ISL $ given by
  $$
  X = \ifundef {varnothing} \emptyset \else \varnothing \fi, \quad Y = \ifundef {varnothing} \emptyset \else \varnothing
\fi \and Z = \big \{(\vr ,1,\vr ): \vr \in E^0\big \},
  $$
  notice that, according to  \cite [Definition 11.4]{\actions }, one has that
  $$
  \ISL ^{X,Y}=\ISL ,
  $$
  and that $Z$ is a cover for $\ISL ^{X,Y}$, as seen in \lcite {\FiniteCover }. By the tightness condition \cite [Definition
11.6]{\actions } we have
  $$
  \bigvee_{z \in Z}\urep (z) \geq \bigwedge_{\vr \in X} \urep (\vr ) \wedge\bigwedge_{y \in Y} \neg{\urep (y)}.
  $$
  As explained in the discussion following \cite [Definition 11.6]{\actions }, the right-hand-side above must be interpreted
as  1 because $X$ and $Y$ are empty.  On the other hand, since the $\urep (z)$ are pairwise orthogonal, the supremum in
the left-hand-side above becomes a sum, so
  $$
  1 = \sum_{z \in Z}\urep (z) = \sum_{\vr \in E^0}\urep (\vr , 1, \vr ) = \tu _1.
  $$

In order to prove that $\tu $ is multiplicative, let $g$ and $h$ be in $G$.  Then
  $$
  \tu_g \tu_h =
  \sum_{\vr ,\vro \in E^0}\urep \big((\vr , g, g^{-1}  \vr )  (\vro , h, h^{-1}  \vro )\big) \$=
  \sum_{\vr \in E^0}\urep \big((\vr , g, g^{-1}  \vr )  (g^{-1}  \vr , h, h^{-1}g^{-1}  \vr )\big) =
  \sum_{\vr \in E^0}\urep \big(\vr , gh, (gh)^{-1}  \vr \big)=
  \tu_{gh}.
  $$

We next claim that  $\tu_g^* =\tu_{g^{-1}}$, for all $g$ in $G$.  To prove it we compute
  $$
  \tu_g^* =
  \sum_{\vr \in E^0}\urep (\vr , g, g^{-1}  \vr )^* =
  \sum_{\vr \in E^0}\urep (g^{-1}  \vr , g^{-1}, \vr ) = \cdots
  $$
  which, upon the change of variables $\vro =g^{-1} \vr $, becomes
  $$
  \cdots = \sum_{\vro \in E^0}\urep (\vro , g^{-1}, gy) = \tu_{g^{-1}}.
  $$

This shows that $\tu $ is a unitary representation, verifying \lcite {\DefineOGE .b}.
Turning now our attention to \lcite {\DefineOGE .c}, let $g \in G$ and $\ed \in E¹$.  Then
  $$
  \tu_g\ts_\ed =
  \sum_{\vr \in E^0}\urep (\vr , g, g^{-1}  \vr )\, \urep \big(\ed , 1, \src (\ed )\big) =
  \urep \big(g\ran (\ed ), g, \ran (\ed )\big)\,\urep \big(\ed , 1, \src (\ed )\big) \$=
  \urep \big(\ran (g\ed )g\ed , \varphi(g,\ed ), \src (\ed )\big) =
  \urep \big(g\ed , \varphi(g,\ed ), \src (\ed )\big) = (\star).
  $$
  On the other hand
  $$
  \ts_{g\ed }\tu_{\varphi(g,\ed )} =
  \urep \big(g\ed ,1,\src (g\ed )\big) \sum_{\vr \in E^0} \urep \big(\vr , \varphi(g,\ed ), \varphi(g,\ed )^{-1}  \vr \big)  \$=
  \urep \big(g\ed ,1,\src (g\ed )\big) \,\urep \big(\src (g\ed ), \varphi(g,\ed ), g^{-1}  \src (g\ed )\big)  \$=
  \urep \big(g\ed ,1,\src (g\ed )\big)\,\urep \big(\src (g\ed ), \varphi(g,\ed ),   \src (\ed )\big)  =
  \urep \big(g\ed , \varphi(g,\ed ), \src (\ed )\big),
  $$
  which coincides with $(\star)$ and hence proves \lcite {\DefineOGE .c}.  We leave the proof of \lcite {\DefineOGE .d} to
the reader after which the universal property of $\OGE $ intervenes to provide us with a *-homomorphism
  $$
  \psi:\OGE \to A
  $$
  sending
  $$
  p_\vr \mapsto \tp_\vr ,\quad \s_\ed \mapsto \ts_\ed \and u_g \mapsto \tu_g.
  $$

Now we must show that
  $$
  \psi\big(\pi(\gamma)\big)=\urep (\gamma) \for \gamma \in \SGE .
  \equationmark PsiPiRho
  $$
  We will first do so for the following
special cases:
  \Zitem $\gamma = (\vr ,1,\vr )$, for $\vr \in E^0$,
  \zitem $\gamma = \big(\ed , 1, \src (\ed )\big)$, for $\ed \in E^1$,
  \zitem $\gamma = (\vr ,g,g^{-1}\vr )$, for $\vr \in E^0$, and $g \in G$.

\bigskip \noindent
In case (i)  we have
  $$
  \psi(\pi(\gamma)) =
  \psi(\pi(\vr ,1,\vr )) =
  \psi(p_\vr ) =
  \tp_\vr =
  \urep (\vr ,1,\vr ) =
  \urep (\gamma).
  $$
  As for (ii),
  $$
  \psi(\pi(\gamma)) =
  \psi\big(\pi\big(\ed , 1, \src (\ed )\big)\big) =
  \psi(\s_\ed ) =
  \ts_\ed =
  \urep \big(\ed , 1, \src (\ed )\big) =
  \urep (\gamma).
  $$
  Under (iii),
  $$
  \psi(\pi(\gamma)) =
  \psi\big(\pi(\vr ,g,g^{-1}\vr ) \big) =
  \psi(p_\vr u_g p_{g^{-1}\vr }) =
  \psi(p_\vr u_g) =
  \tp_\vr \tu_g \$=
  \urep (\vr ,1,\vr ) \sum_{\vro \in E^0}\urep (\vro , g, g^{-1} \vro ) =
  \sum_{\vro \in E^0}\urep \big((\vr ,1,\vr ) (\vro , g, g^{-1} \vro )\big) =
  \urep (\vr , g, g^{-1} \vr ) =
  \urep (\gamma).
  $$

In order to prove \lcite {\PsiPiRho }, it is now clearly enough to check that the *-sub-semigroup of $\SGE $ generated
by the elements mentioned in (i--iii), above, coincides with $\SGE $.

Denoting this *-sub-semigroup by $\T $, we will first show that $\big(\alpha,1,\src (\alpha)\big)$ is in $\T $, for
every $\alpha \in E^*$.  This is evident for $|\alpha|\leq1$, so we suppose that $\alpha=\alpha'\alpha''$, with $\alpha'
\in E^1$, and $\ran (\alpha'')=\src (\alpha')$.  We then have by induction that
  $$
  \T \ni
  \big(\alpha',1,\src (\alpha')\big)   \big(\alpha'',1,\src (\alpha'')\big) =
  \big(\alpha\alpha'',1,\src (\alpha'')\big) =
  \big(\alpha,1,\src (\alpha)\big).
  $$


Considering a general element $(\alpha,g,\beta) \in \SGE $,
  let $\vr =\src (\alpha)$, so that $g^{-1}\vr =\src (\beta)$, and notice that
  $$
  \T \ni
  \big(\alpha,1,\src (\alpha)\big) (\vr , g, g^{-1}  \vr ) \big(\beta,1,\src (\beta)\big)^* \$=
  \big(\alpha,1,\src (\alpha)\big) \big(\src (\alpha), g, \src (\beta)\big) \big(\src (\beta),1,\beta\big) =
  (\alpha,g,\beta),
  $$
  which proves that $\T =\SGE $, and hence that \lcite {\PsiPiRho } holds.

To conclude we observe that the uniqueness of $\psi$ follows  from the fact that $\OGE $ is generated by the $p_\vr $, the
$\s_\ed $, and the $u_g$.
  \endProof

Given an inverse semigroup $\S $ with zero, recall from \cite [Theorem 13.3]{\actions } that $\G \tight (\S )$ (denoted simply as
$\G \tight $ in \cite {\actions }) is the groupoid of germs for the natural action of $\S $ on the space of tight filters over
its idempotent semi-lattice.  Moreover the C*-algebra of $\G \tight (\S )$ is universal for tight representations of $\S $.

\state Corollary Under the assumptions of \lcite {\StandingHyp } one has that $\OGE $ is isomorphic to the C*-algebra of
the groupoid $\Gpd $.

\Proof Follows from \cite [Theorem 13.3]{\actions } and the uniqueness of universal C*-algebras. \endProof

We should notice that our requirement that $G$ be countable in \lcite{\StandingHyp} is only used in the above proof,
where the application of \cite [Theorem 13.3]{\actions } depends on the countability of $\SGE$.

We conclude by analyzing  the question of Hausdorffness for $\Gpd $.

\state Proposition When $(\Data )$ is {\essfree }, one has that  $\Gpd $ is a Hausdorff groupoid.

\Proof By \lcite {\EssFreeEUnitary } we have that $\SGE $ is E*-unitary.  The result then follows from \cite [Propositions
6.2 \& 6.4]{\actions }.
  \endProof

\section Corona Groups

It is our next goal to give a concrete description of $\Gpd $, similar to the description given to the groupoid associated
to a row-finite graph in  \cite [Definition 2.3]{\KPRR}.  The crucial ingredient there is the notion of \"{tail equivalence with
lag}.  In this section we will construct a group where  our generalized \"{lag} function will take values.

Let $G$ be a group.  Within the infinite cartesian product\fn{For the purposes of this construction we adopt the
convention that ${\bf N}=\{1,2,3,\ldots\}$.}
  $$
  G^\infty = \prod_{n \in {\bf N}} G
  $$
  consider the infinite direct sum
  $$
  G^{(\infty)} = \bigoplus_{n \in {\bf N}} G
  $$
  formed by the elements  $g=(g_n)_{n \in {\bf N}} \in G^\infty$ which are eventually trivial, that is,
for which there exists $n_0$ such that $g_n=1$, for all $n\geq n_0$.
  It is clear that $G^{(\infty)}$ is a normal subgroup of $G^\infty$.

\definition Given a group $G$, the \"{corona} of $G$ is the quotient group
  $$
  \corona = G^\infty/G^{(\infty)}.
  $$

  Consider the \"{left} and \"{right shift} endomorphisms of $G^\infty$
  $$
  \lambda,\rho : G^\infty\to G^\infty
  $$
  given for every $\g =(\g_n)_{n \in {\bf N}} \in G^\infty$, by
  $$
  \lambda(\g )_n = \g_{n+1} \for n \in {\bf N},
  $$
  and
  $$
  \rho(\g )_n = \left \{ \matrix {1, & \hbox { if } n=0, \cr \g_{n-1}, & \hbox { if } n\geq1.} \right .
  $$
  It is readily seen that $G^{(\infty)}$ is invariant under both $\lambda$ and $\rho$, so these pass to the quotient providing
endomorphisms
  $$
  \q \lambda, \q \rho : \corona \to \corona .
  \equationmark DefineQuotientAuto
  $$

For every  $\g =(\g_n)_{n \in {\bf N}} \in G^\infty$, we have that
  $$
  \lambda(\rho(\g ))= \g
  \and
  \rho(\lambda(\g )) = (1,\g _2,\g _3,\ldots) \equiv  \g ,
  \equationmark ShiftIsAuto
  $$
  where we use ``$\equiv$'' to refer to the equivalence relation determined by the normal subgroup $G^{(\infty)}$.
  Therefore both $\q \lambda\q \rho$ and $\q \rho\q \lambda$ coincide with the identity, and hence $\q \lambda$ and $\q \rho$ are
each other's inverse.  In particular, they are both automorphisms of $\corona $.

Iterating $\q \rho$ therefore gives an action of ${\bf Z}$ on $\corona $.

\definition Given any countable discrete group $G$, the \"{lag group} associated to $G$ is the  semi-direct product group
  $$
  \corona \ifundef {rtimes} \times \else \rtimes \fi_{\q \rho} {\bf Z}.
  $$

The reason we call this the ``lag group'' is that it will play a very important role in the next section, as the co-domain
for our \"{lag} function.

\section The tight groupoid of $\SGE $

We would now like to give a detailed description of the groupoid $\Gpd $.  As already mentioned this is the
groupoid of germs for the natural action of $\SGE $ on the space of tight filters over the idempotent semi-lattice
$\ISL $ of $\SGE $.  See \cite [Section 4]{\actions } for more details.

Given an infinite word
  $$
  \xi=\xi_1\xi_2\ldots \in E^\infty,
  $$
  and an integer $n\geq0$, denote by $\trunc \xi n$ the finite word of length $n$ given by
  $$
  \trunc \xi n = \left \{\matrix {\xi_1\xi_2\ldots\xi_n, & \hbox { if } n\geq1, \cr \cr
  \ran (\xi_1),  & \hbox { if } n=0.}\right .
  $$

\state Proposition \label ActionOnInfWords
  There is a unique action
  $$
  (g,\xi) \in G\times E^\infty \mapsto g\xi \in E^\infty
  $$
  of $G$ on $E^\infty$ such that,
  $$
  \trunc {(g\xi)}n = g (\trunc \xi n),
  $$
  for every $g \in G$, $\xi \in E^\infty$, and $n \in {\bf N}$.

\Proof Left to the reader. \endProof

Recall from \lcite {\DefineEAlpha } that, for any finite word $\alpha \in E^*$, we denote by $\eproj_\alpha$ the
idempotent element $(\alpha,1,\alpha)$ in $\ISL $.  Thus, given an infinite word $\xi \in E^\infty$, we may look at the
subset
  $$
  {\cal F}_\xi = \{\eproj_{\xi_n}: n \in {\bf N}\} \subseteq \ISL ,
  $$
  which turns out to be an ultra-filter \cite [Definition]{\actions } over $\ISL $.  Denoting the set of all ultra-filters
over $\ISL $ by $\widehat \ISL_\infty$, as in \cite [Definition 12.8]{\actions }, one may also show
  \cite [Proposition 19.11]{\actions }
  that the  correspondence
  $$
  \xi \in E^\infty \mapsto {\cal F}_\xi \in \widehat \ISL_\infty
  $$
  is bijective, and we will use it to identify
  $
  E^\infty
  $
  and
  $
  \widehat \ISL_\infty.
  $
  Furthermore, this correspondence may be proven to be a homeomorphism if $E^\infty$ is equipped with the product topology.

Since $E$ is finite, $E^\infty$ is compact by Tychonov's Theorem, and consequently so is $\widehat \ISL_\infty$.  Being the closure
of $\widehat \ISL_\infty$ within $\widehat \ISL $ \cite [Theorem 12.9]{\actions }, the space $\widehat \ISL \tight $ formed by the tight filters
therefore necessarily coincides with $\widehat \ISL_\infty$.

Identifying $\widehat \ISL \tight $ with $E^\infty$, as above, we may transfer the canonical action of $\SGE $ from the former
to the latter resulting in the following: to each element
  $(\alpha,g,\beta) \in \SGE $,
  we associate the partial homeomorphism of $E^\infty$ whose domain is the \"{cylinder}
  $$
  \cyl \beta := \{\eta \in E^\infty: \eta =\beta\xi, \hbox { for some } \xi \in E^\infty\},
  $$
  and which sends each $\eta =\beta\xi \in \cyl \beta$ to $\alpha g\xi $, where the meaning of ``$g\xi $''  is as in \lcite {\ActionOnInfWords }.

  As before we will not use any special symbol to indicate this action, using module notation instead:
  $$
  (\alpha,g,\beta)\eta  = \alpha g\xi .
  $$

Before we proceed let us at least check that $\alpha g\xi $ is in fact an element of $E^\infty$, which is to say that
  $
  \src (\alpha) = \ran (g\xi).
  $
  Firstly, for every element $(\alpha,g,\beta) \in \SGE $, we have that $\src (\alpha)=g\src (\beta)$.  Secondly, if
$\eta =\beta\xi \in E^\infty$, then $\src (\beta) = \ran (\xi)$.  Therefore
  $$
  \ran (g\xi) = g\ran (\xi) = g\src (\beta) = \src (\alpha).
  $$

  This leads to a first, more or less concrete description of $\Gpd $, namely the groupoid of germs for the
above action of $\SGE $ on $E^\infty$.  Our aim is nevertheless a much more precise description of it.

Recall from \cite [Definition 4.6]{\actions } that the germ of an element $s \in \SGE $ at a point $\xi$ in the domain of $s$  is denoted by
  $[s,\xi]$.  If $s=(\alpha,g,\beta)$, this would lead to the somewhat awkward notation
  $[(\alpha,g,\beta),\xi]$, which from now on will be simply written as
  $$
  \germ \alpha g\beta\xi.
  $$

Thus the groupoid $\Gpd $, consisting of all germs for the action of $\SGE $ on $E^\infty$, is given by
  $$
  \Gpd =
  \Big \{\germ \alpha g\beta\xi: (\alpha,g,\beta) \in \SGE ,\ \xi \in \cyl \beta \Big \}.
  \equationmark FirstModel
  $$

Let us now prove a criterion for equality of germs.

\state Proposition \label EqualGerms
  Suppose that $(\Data )$ is {\essfree } and let us be given elements $(\alpha_1,g_1,\beta_1)$ and $(\alpha_2,g_2,\beta_2)$ in
$\SGE $, with $|\beta_1|\leq|\beta_2|$, as well as infinite paths $\eta _1$ in $\cyl {\beta_1}$, and $\eta _2$ in $\cyl {\beta_2}$.  Then
  $$
  \germ {\alpha_1}{g_1}{\beta_1}{\eta _1}  =  \germ {\alpha_2}{g_2}{\beta_2}{\eta _2}
  $$
  if and only if there is a finite path $\gamma \in E^*$ and an infinite path $\xi \in E^\infty$, such that
  \Zitem $\alpha_2 = \alpha_1 g_1\gamma,$
  \zitem $g_2 = \varphi(g_1,\gamma),$
  \zitem $\beta_2 = \beta_1\gamma,$
  \zitem $\eta _1=\eta _2=\beta_1\gamma\xi .$

\Proof Assuming that the germs are equal, we have by \cite [Definition 4.6]{\actions } that
  $$
  \eta _1=\eta _2 =: \eta,
  $$
  and there is an
idempotent $(\delta,1,\delta) \in \ISL $, such that $\eta \in \cyl \delta$, and
  $$
  (\alpha_1,g_1,\beta_1) (\delta,1,\delta) = (\alpha_2,g_2,\beta_2) (\delta,1,\delta).
  \equationmark SameGerm
  $$

  It follows that $\eta =\delta\zeta$, for some $\zeta \in E^\infty$.  Upon replacing $\delta$ by a longer prefix of $\eta $, we may assume that $|\delta|$ is as
large as we want.  Furthermore the element of $\SGE $ represented by the two sides of
the equation displayed above is evidently nonzero because the partial homeomorphism associated to it under our
action has $\eta $ in its domain.  So, focusing on \lcite {\inverseSemigroup }, we see that $\beta_1$ and $\delta$ are comparable, and
so are $\beta_2$ and $\delta$.

Assuming that $|\delta|$ exceeds both $|\beta_1|$ and $|\beta_2|$, we may then write $\delta = \beta_1\varepsilon_1 = \beta_2\varepsilon_2$, for suitable $\varepsilon_1$ and $\varepsilon_2$ in
$E^*$.  But since $|\beta_1|\leq|\beta_2|$, this in turn implies that $\beta_2=\beta_1\gamma$, for some $\gamma \in E^*$, hence proving (iii).
  Therefore
  $
  \delta=\beta_1\gamma\varepsilon_2,
  $
  so
  $$
  \eta  = \delta\zeta =  \beta_1\gamma\varepsilon_2\zeta,
  $$
  and (iv) follows once we choose $\xi  = \varepsilon_2\zeta$.  Moreover, equation \lcite {\SameGerm } reads
  $$
  (\alpha_1,g_1,\beta_1) (\beta_1\gamma\varepsilon_2,1,\beta_1\gamma\varepsilon_2) = (\alpha_2,g_2,\beta_1\gamma) (\beta_1\gamma\varepsilon_2,1,\beta_1\gamma\varepsilon_2).
  $$

  Computing the products according to \lcite {\inverseSemigroup }, we get
  $$
  \big(\alpha_1g_1(\gamma\varepsilon_2),\ \varphi(g_1,\gamma\varepsilon_2),\ \beta_1\gamma\varepsilon_2\big) = \big(\alpha_2g_2\varepsilon_2,\ \varphi(g_2,\varepsilon_2),\ \beta_1\gamma\varepsilon_2\big),
  $$
  from where we obtain
  $$
  \alpha_2g_2\varepsilon_2 =  \alpha_1g_1(\gamma\varepsilon_2) = \alpha_1(g_1\gamma)\varphi(g_1,\gamma)\varepsilon_2,
  \equationmark EqOne
  $$
  and
  $$
  \varphi(g_2,\varepsilon_2) = \varphi(g_1,\gamma\varepsilon_2)  =  \varphi\big(\varphi(g_1,\gamma),\varepsilon_2\big).
  \equationmark EqTwo
  $$

  Since $|g_2\varepsilon_2| = |\varepsilon_2| = |\varphi(g_1,\gamma)\varepsilon_2|$, we deduce from \lcite {\EqOne } that
  $$
  g_2\varepsilon_2 = \varphi(g_1,\gamma)\varepsilon_2,
  \equationmark DeduceOne
  $$
  and hence also that
  $$
  \alpha_2 = \alpha_1g_1\gamma,
  $$
  proving (i).
  Defining $g=g_2^{-1}\varphi(g_1,\gamma)$, we claim that
  $$
  g\varepsilon_2=\varepsilon_2\and \varphi(g,\varepsilon_2) = 1.
  $$

In view of \lcite {\EqTwo } and \lcite {\DeduceOne }, point (ii) follows from \lcite {\VarStar }.

Conversely, assume (i--iv) and let us prove equality of the above germs.  Setting $\delta = \beta_1\gamma$, we have by (iv) that
  $$
  \eta := \eta _1 = \eta _2 \in \cyl \delta,
  $$
  so it suffices to verify \lcite {\SameGerm }, which the reader could do without any difficulty.
\endProof

The above result then says that the typical situation in which an  equality of germs takes place is
  $$
  \germ {\alpha}{g}{\beta}{\beta\gamma\xi }  =  \germ {\alpha g\gamma}{\varphi(g,\gamma)}{\beta\gamma}{\beta\gamma\xi }.
  $$

Our next two results are designed to offer convenient representatives of germs.

\state Proposition \label AnyLength
  Given any germ $u$, there exists an integer $n_0$, such that for every $n\geq n_0$,
  \Zitem there is a representation of $u$ of the form $u=\germ {\alpha_1} {g_1} {\beta_1} {\beta_1\xi_1}$, with $|\alpha_1|=n$.
  \zitem there is a representation of $u$ of the form $u=\germ {\alpha_2} {g_2} {\beta_2} {\beta_2\xi_2}$, with $|\beta_2|=n$.

\Proof
  Write $u=\germ \alpha g \beta \eta$, and choose any  $n_0\geq\max \{|\alpha|,|\beta|\}$.  Then, for every $n\geq n_0$ we may write $\eta=\beta\gamma\xi$, with $\gamma \in
E^*$, $\xi \in E^\infty$, and
such that $|\gamma|=n-|\alpha|$ (resp.~$|\gamma|=n-|\beta|$).  Therefore
  $$
  u = \germ \alpha g \beta {\beta\gamma\xi} = \germ {\alpha g\gamma} {\varphi(g,\gamma)} {\beta\gamma} {\beta\gamma\xi},
  $$
  and we have
  $
  |\alpha g\gamma| = |\alpha| + |g\gamma| = |\alpha| + |\gamma| = n
  $
  (resp.~
  $
  |\beta\gamma| = |\beta| + |\gamma| = n
  $).
  \endProof

\state Corollary \label MultPairs
  Given $u_1$ and $u_2$ in $\Gpd $, such that $(u_1,u_2) \in \Gtwo $ (that is, such that the multiplication $u_1u_2$ is allowed
or, equivalently, such that $\src(u_1)=\ran(u_2)$), there are representations of $u_1$
and $u_2$ of the form
  $$
  u_1 = \germ {\alpha_1} {g_1} {\alpha_2} {\alpha_2g_2\xi } \and
  u_2 = \germ {\alpha_2} {g_2} {\beta } {\beta \xi },
  $$
  and in this case
  $$
  u_1u_2 =   \germ {\alpha_1} {g_1g_2} {\beta } {\beta \xi }.
  $$

\Proof
  Using \lcite {\AnyLength }, write
  $$
  u_i =
  \germ {\alpha_i} {g_i} {\beta_i} {\beta_i\xi_i},
  $$
  with $|\beta_1|=|\alpha_2|$.
  By virtue of $(u_1,u_2)$ lying in $\Gtwo $, we have that
  $$
  \beta_1\xi_1 = (\alpha_2,g_2,\beta_2)(\beta_2\xi_2) = \alpha_2g_2\xi_2,
  $$
  so in fact $\beta_1 = \alpha_2$, and $\xi_1 = g_2\xi_2$.  Then
  $$
  u_1 =
  \germ {\alpha_1} {g_1} {\beta_1} {\beta_1\xi_1} =
  \germ {\alpha_1} {g_1} {\alpha_2} {\alpha_2g_2\xi_2},
  $$
  and it suffices to put $\xi =\xi_2$, and $\beta = \beta_2$.

With respect to the last assertion we have that $u_1u_2 = [s;\beta\xi]$, where $s$ is the element of $\SGE $ given by
  $$
  s =
  (\alpha_1, g_1, \alpha_2) (\alpha_2, g_2,\beta) \={\EasyMult }
  (\alpha_1, g_1g_2, \beta),
  $$
  concluding the proof.
  \endProof

Having extended the action of $G$ to the set of infinite words in \lcite {\ActionOnInfWords }, one may ask whether it is
possible to do the same for the cocycle $\varphi$.  The following is an attempt at this which however produces a map taking
values in the infinite product $G^\infty$, rather than in $G$.

  \definition We will denote by $\Phi$, the map
  $$
  \Phi: G\times E^\infty \to G^\infty
  $$
  defined by the rule
  $$
  \Phi(g,\xi)_n=\varphi(g,\trunc \xi{n-1}),
  $$
  for $g\in G$, $\xi\in E^\infty$, and $n\geq1$.

We wish to view $\Phi$ as some sort of cocycle but, unfortunately, property \lcite {\Equacoes .\LastItem } does not hold quite
as stated.  On the fortunate side, a suitable modification of $\Phi$, involving the left shift endomorphism $\lambda$ of $G^\infty$,
works nicely:

\state Proposition \label LastWithShift
  Let $\alpha$ be a finite word in $E^*$ and let $\xi$ be an infinite word in $E^\infty$ such that $\src (\alpha)=\ran (\xi)$.  Then, for
every $g$ in $G$, one has that
  $$
  \Phi\big(\varphi(g,\alpha),\xi\big) =
  \lambda^{|\alpha|}\big(\Phi\big(g,\alpha\xi)\big)
  $$

\Proof For all $n\geq1$, we have
  $$
  \Phi\big(\varphi(g,\alpha),\xi\big)_n =
  \varphi\big(\varphi(g,\alpha),\trunc \xi{n-1}\big) =
  \varphi\big(g,\alpha(\trunc \xi{n-1})\big) \$=
  \varphi\big(g,\trunc {(\alpha\xi)}{n-1+|\alpha|}\big) =
  \lambda^{|\alpha|}\big(\Phi(g,\alpha\xi)\big)_n.
  \endProof

Another reason to think of  $\Phi$ as a cocycle  is the following version of the cocycle identity \lcite {\Equacoes .b}:

\state Proposition \label CocycleIdForCapPhi
  For every $\xi \in E^\infty$, and every $g,h \in G$, we have that
  $$
  \Phi(gh,\xi) =  \Phi\big(g,h\xi\big)\Phi(h,\xi).
  $$

\Proof We have for all $n \in {\bf N}$, that
  $$
  \Phi(gh,\xi)_n =
  \varphi(gh,\xi|_{n-1}) \={\Equacoes .b}
  \varphi\big(g,h(\xi|_{n-1})\big)\varphi(h,\xi|_{n-1}) \={\ActionOnInfWords } $$$$=
  \varphi\big(g,(h\xi)|_{n-1}\big)\varphi(h,\xi|_{n-1}) =
  \Phi\big(g,h\xi\big)_n\Phi(h,\xi)_n.
  \endProof

The following elementary fact might perhaps justify the choice of ``$n-1$'' in the definition of $\Phi$.

\state Proposition \label ActionCocycle
  Given $g \in G$, and $\xi \in E^\infty$, one has that
  $$
  (g\xi)_n = \Phi(g,\xi)_n\,\xi_n.
  $$

  \Proof
  By \lcite {\ActionOnInfWords } we have that  $(g\xi)|_n = g(\xi|_n)$, so the  $n^{th}$ letter of $g\xi$ is also the $n^{th}$
letter of $g(\xi|_n)$.  In addition we have that
  $$
  g(\xi|_n) =
  g(\xi|_{n-1}\xi_n) \={\Equacoes .\MainConcat }
  g(\xi|_{n-1})\varphi(g,\xi|_{n-1})\xi_n,
  $$
  so
  $$
  (g\xi)_n = \varphi(g,\xi|_{n-1})\xi_n = \Phi(g,\xi)_n\,\xi_n.
  \endProof

We now wish to define a  homomorphism (also called a one-cocycle) from  $\Gpd $ to the lag group   $\corona \ifundef {rtimes} \times \else \rtimes \fi_\rho {\bf Z}$, by means of the
rule
  $$
  \germ \alpha g\beta{\beta\xi} \mapsto \Big (\rho^{|\alpha|}\big(\Phi(g,\xi)\big),|\alpha|-|\beta|\Big ).
  $$

  As it is often the case for maps defined on groupoid of germs, the above tentative definition uses a representative of
the germ, so some work is necessary to prove that the definition does not depend on the choice of the representative.
The technical part of this task is the content of our next result.

\state Lemma \label LagWellDefined
  Suppose that $(\Data )$ is {\essfree }.  For each $i=1,2$, let us be given $(\alpha_i,g_i,\beta_i)$ in $\SGE $, as well as $\eta
_i=\beta_i\xi_i \in \cyl {\beta_i}$.  If
  $$
  \germ {\alpha_1}{g_1}{\beta_1}{\eta _1}  =  \germ {\alpha_2}{g_2}{\beta_2}{\eta _2},
  $$
  then
  $$
  \rho^{|\alpha_1|}\big(\Phi(g_1,\xi_1)\big) \equiv \rho^{|\alpha_2|}\big(\Phi(g_2,\xi_2)\big)
  $$
  modulo  $G^{(\infty)}$.

\Proof
  Assuming without loss of generality that $|\beta_1|\leq|\beta_2|$, we may use \lcite {\EqualGerms } to write
  $$
  \alpha_2 = \alpha_1 g_1\gamma, \quad g_2 = \varphi(g_1,\gamma), \quad \beta_2 = \beta_1\gamma \and \eta _1=\eta_2=\beta_1\gamma\xi,
  $$
  for suitable  $\gamma \in E^*$ and $\xi \in E^\infty$.  Then necessarily $\xi_1 = \gamma\xi$, and $\xi_2=\xi$, and
  $$
  \rho^{|\alpha_2|}\big(\Phi(g_2,\xi_2)\big) =
  \rho^{|\alpha_1|+|\gamma|}\Big ( \Phi\big(\varphi(g_1,\gamma_1),\xi\big)\Big ) \={\LastWithShift } $$ $$=
  \rho^{|\alpha_1|}\rho^{|\gamma|}\lambda^{|\gamma|}\big(\Phi\big(g_1,\gamma\xi)\big)  \explain \ShiftIsAuto \equiv
  \rho^{|\alpha_1|}\big(\Phi\big(g_1,\xi_1)\big).
  \endProof

\fix Due to our reliance on
\lcite {\EqualGerms } and
\lcite {\LagWellDefined },
from now on and until the end of this section
we will assume, in addition to  \lcite {\StandingHyp }, that $(\Data )$ is {\essfree }.

\bigskip
If $\g $ is in $G^\infty$, we will denote by $\q \g $ its class in the quotient group $\corona $.  Likewise we will denote by
$\q \Phi$ the composition of $\Phi$ with the quotient map

  \beginpicture
  \setcoordinatesystem units <0.0040truecm, -0.0040truecm> point at 3000 0
  \setplotarea x from -1500 to 1000, y from -200 to 500
  \put {$G\times E^\infty \longrightarrow \ G^\infty \longrightarrow \corona $} at 510 000
  \put {$\Phi$} at 440 -90
  \setquadratic
  \plot 100 120 525 250 950 120 /
  \arrow <0.11cm> [0.3,1.2] from 950 120 to 960 113
  \put {$\q \Phi$} at 525 350
  \endpicture

  \bigskip \noindent from $G^\infty$ to $\corona $.  It then follows from the above Lemma that the correspondence
  $$
  \germ \alpha g\beta{\beta\xi} \in \Gpd \mapsto  \q \rho^{|\alpha|}\big(\q \Phi(g,\xi)\big) \in \corona
  $$
  is a well defined map.  This is an important part of the one-cocycle we are about to introduce.

\state Proposition \label DefineLag
  The correspondence
  $$
  \lag :   \germ \alpha g\beta{\beta\xi} \mapsto \Big (\q \rho^{|\alpha|}\big(\q \Phi(g,\xi)\big),|\alpha|-|\beta|\Big )
  $$
  gives a well defined map
  $$
  \lag : \Gpd \to \corona \ifundef {rtimes} \times \else \rtimes \fi_\rho {\bf Z},
  $$
  which is moreover a one-cocycle.  From now on $\lag$ will be called the \"{lag function}.

  \Proof
  By the discussion above we have that the first coordinate of the above pair is well defined.  On the other hand, in
the context of \lcite {\EqualGerms } one easily sees that $|\alpha_1|-|\beta_1| = |\alpha_2|-|\beta_2|$, so the second coordinate is also well
defined.

In order to show that $\lag $ is multiplicative, pick $(u_1,u_2) \in \Gtwo $.  We may then use \lcite {\MultPairs } to write
  $$
  u_1 = \germ {\alpha_1} {g_1} {\alpha_2} {\alpha_2g_2\xi }
  \and
  u_2 = \germ {\alpha_2} {g_2} {\beta } {\beta \xi }.
  $$
  So
  $$
  \lag (u_1)\lag (u_2) =
  \Big (\rho^{|\alpha_1|}\big(\Phi(g_1,g_2\xi)\big),\ |\alpha_1|-|\alpha_2|\Big )
  \Big (\rho^{|\alpha_2|}\big( \Phi(g_2,\xi)\big),\ |\alpha_2|-|\beta|\Big )
  \$=
  \Big (\rho^{|\alpha_1|}\big(\Phi(g_1,g_2\xi)\big) \ \rho^{|\alpha_1|}\big( \Phi(g_2,\xi)\big),\ \ |\alpha_1|-|\alpha_2|+|\alpha_2|-|\beta|\Big )  \$=
  \Big (\rho^{|\alpha_1|}\Big (\Phi(g_1,g_2\xi) \Phi(g_2,\xi)\Big ),\ |\alpha_1|-|\beta|\Big )  \={\CocycleIdForCapPhi }
  \Big (\rho^{|\alpha_1|}\big(\Phi(g_1g_2,\xi)\big),\ |\alpha_1|-|\beta|\Big ) \$=
  \lag \big(\germ {\alpha_1} {g_1g_2} {\beta } {\beta \xi }\big) \=\MultPairs
  \lag (u_1u_2).
  \endProof

The main relevance of this one-cocycle is that it essentially describes the elements of $\Gpd $, as we would like to show
now.

\state Proposition \label FInjective
  Given
  $
  u_1,u_2 \in \Gpd ,
  $
  one has that
  $$
  \def \ker {\kern 8pt}
  \def \and {\ker \wedge\ker }
  \src (u_1)=\src (u_2) \and \ran (u_1)=\ran (u_2) \and \lag (u_1)=\lag (u_2) \imply u_1=u_2.
  $$

\Proof Using \lcite {\AnyLength }, write $u_i = \germ {\alpha_i}{g_i}{\beta_i}{\beta_i\xi_i}$, for $i=1,2$, with $|\beta_1|=|\beta_2|$.  Since
  $$
  \beta_1\xi_1 = \src (u_1) = \src (u_2) =  \beta_2\xi_2,
  $$
  we conclude that
  $
  \beta_1=\beta_2,
  $
  and
  $$
  \xi_1=\xi_2=:\xi.
  $$

By focusing on the second coordinate of $\lag (u_i)$, we see that $|\alpha_1|-|\beta_1|=|\alpha_2|-|\beta_2|$, and hence $|\alpha_1|=|\alpha_2|$.  Moreover,
since
  $$
  \alpha_1g_1\xi =  \alpha_1g_1\xi_1 = \ran (u_1) =   \ran (u_2) =  \alpha_2g_2\xi_2 = \alpha_2g_2\xi,
  $$
  we see that $\alpha_1=\alpha_2$, and
  $$
  g_1\xi  = g_2\xi.
  \equationmark EqualGXi
  $$

The fact that $\lag (u_1) = \lag (u_2)$ also implies that
  $$
  \q \rho^{|\alpha_1|}\big(\q \Phi(g_1,\xi)\big) =
  \q \rho^{|\alpha_2|}\big(\q \Phi(g_2,\xi)\big),
  $$
  and since $\alpha_1=\alpha_2$, we conclude that
  $
  \q \Phi(g_1,\xi) = \q \Phi(g_2,\xi),
  $
  and hence that there exists an integer $n_0$ such that
  $$
  \varphi(g_1,\xi|_n) = \varphi(g_2,\xi|_n)
  \for n\geq n_0.
  $$
  By \lcite {\EqualGXi } we also have that
  $
  g_1(\xi|_n)  = g_2(\xi|_n),
  $
  so \lcite {\VarStar } gives $g_1=g_2$, whence $u_1=u_2$.
  \endProof

As a consequence of the above result we see that the map
  $$
  F:\Gpd \to E^\infty\times (\corona \ifundef {rtimes} \times \else \rtimes \fi_{\q \rho} {\bf Z}) \times E^\infty
  $$
  defined by the rule
  $$
  F(u) = \big(\ran (u),\lag (u),\src (u)\big),
  \equationmark DefineF
  $$
  is one-to-one.

Observe that the co-domain of $F$ has a natural groupoid structure, being  the cartesian product of the lag group $\corona
\ifundef {rtimes} \times \else \rtimes \fi_{\q \rho} {\bf Z}$ by the graph of the transitive equivalence relation on $E^\infty$.

Putting together \lcite {\DefineLag } and \lcite {\FInjective } we may now easily prove:

\state Corollary $F$ is a groupoid homomorphism (functor), hence establishing  an isomorphism from $\Gpd $ to the range of $F$.

The range of $F$ is then the concrete model of $\Gpd $ we are after.  But, before giving a detailed description of it, let us make a remark concerning notation: since the co-domain of $F$ is a mixture of cartesian and semi-direct
products, the standard notation for its elements would be something like $\big(\eta,(u,p),\zeta\big)$, for $\eta,\zeta \in E^\infty$, $u \in \corona $,
and $p \in {\bf Z}$.  As part of our effort to avoid heavy notation we will instead denote such an element by
  $$
  \big(\eta;u,p;\zeta\big).
  $$

\state Proposition
  The range of $F$ is precisely the subset of $E^\infty\times (\corona \ifundef {rtimes} \times \else \rtimes \fi_{\q \rho} {\bf Z}) \times E^\infty$, formed by the elements
  $
  (\eta;\q \g ,p-q;\zeta),
  $
  where $\eta,\zeta \in E^\infty$, $\g \in G^\infty$, and $p,q \in {\bf N}$, are such that, for all $n\geq1$,
  \Zitem $\g_{n+p+ 1} = \varphi(\g_{n+p},\zeta_{n+q})$,
  \zitem $\eta_{n+p} = \g_{n+p}\zeta_{n+q}$.

\Proof
  Pick a general element $\germ \alpha g\beta{\beta\xi} \in \Gpd $ and, recalling that
  $$
  F(\germ \alpha g\beta{\beta\xi}) = \Big (\alpha g\xi;\ \q \rho^{|\alpha|}\big(\q \Phi(g,\xi)\big),\ |\alpha|-|\beta|;\ \beta\xi\Big ),
  \equationmark GeneralRangeF
  $$
  let
  $
  \eta=\alpha g\xi, \ \g = \rho^{|\alpha|}\big(\Phi(g,\xi)\big), \ p=|\alpha|,\  q=|\beta|, \hbox { and } \zeta=\beta\xi,
  $
  so that the element depicted in \lcite {\GeneralRangeF } becomes   $(\eta;\q \g ,p-q;\zeta)$, and we must now verify (i) and
(ii).
   For all $n\geq1$, one has that
  $$
  \g_{n+|\alpha|} = \Phi(g,\xi)_n = \varphi(g,\xi|_{n-1}),
  $$
  so
  $$
  \eta_{n+p} =
  (\alpha g\xi)_{n+|\alpha|} =
  (g\xi)_n \={\ActionCocycle }
  \varphi(g,\xi|_{n-1})\xi_n =
  \g_{n+|\alpha|}(\beta\xi)_{n+|\beta|} =
  \g_{n+p}\zeta_{n+q},
  $$
  proving (ii).  Also,
  $$
  \g_{n+p+ 1} =
  \g_{n+|\alpha|+ 1} =
  \varphi\big(g,\xi|_{n}\big) =
  \varphi(g,\xi|_{n-1}\xi_{n}) =
  \varphi\big(\varphi(g,\xi|_{n-1}),\xi_{n}\big) \$=
  \varphi\big(\g_{n+|\alpha|},(\beta\xi)_{n+|\beta|}\big) =
  \varphi\big(\g_{n+p},\zeta_{n+q}\big),
  $$
  proving (i) and hence showing that the range of $F$ is a subset of the set described in the statement.

Conversely, pick  $\eta,\zeta \in E^\infty$, $\g \in G^\infty$, and $p,q \in {\bf N}$ satisfying (i) and (ii), and let us  show that the element
  $(\eta;\q \g ,p-q;\zeta)$ lies in the range of $F$.  Let
  $$
  g  = \g_{p+1},\quad
  \alpha = \eta|_p \and \beta=\zeta|_q,
  $$
  so $\zeta=\beta\xi$ for a unique $\xi \in E^\infty$.   We then claim that $\germ \alpha g\beta{\beta\xi}$ lies in  $\Gpd $.  In order to see this notice that
  $$
  g\src (\beta) =
  g\src (\zeta_q) =
  g\ran (\zeta_{q+1}) =
  \ran (\g_{p+1}\zeta_{q+1}) \={ii}
  \ran (\eta_{p+1}) =
  \src (\eta_p) =
  \src (\alpha),
  $$
  so $(\alpha,g,\beta) \in \SGE $, and therefore $\germ \alpha g\beta{\beta\xi}$ is indeed a member of   $\Gpd $.
  The proof will then be concluded once we show that
  $$
  F(\germ \alpha g\beta{\beta\xi}) =  (\eta;\q \g ,p-q;\zeta),
  $$
  which in turn is equivalent to showing that
  \iaitem
  \aitem $\alpha g\xi = \eta,$
  \aitem $\q \rho^{|\alpha|}\big(\q \Phi(g,\xi)\big) = \q \g ,$
  \aitem $|\alpha|-|\beta| = p-q,$
  \aitem $\beta\xi = \zeta$.

\medskip Before proving these points we will show  that
  $$
  \varphi(\g_{p+1},\xi|_n) = \g_{n+p+1} \for n\geq0.
  \eqno {(\dagger)}
  $$
  This is obvious for $n=0$.  Assuming that $n\geq1$ and using induction, we have
  $$
  \varphi(\g_{p+1},\xi|_n) =
  \varphi(\g_{p+1},\xi|_{n-1}\xi_n) =
  \varphi\big(\varphi(\g_{p+1},\xi|_{n-1}),\xi_n\big) \$=
  \varphi(\g_{n+p},\zeta_{n+q}) \={i}   \g_{n+p+1},
  $$
  verifying $(\dagger)$.

Addressing (a) we have to prove that $(\alpha g\xi)_k = \eta_k$, for all $k\geq1$, but given that $\alpha$ is defined to be $\eta|_p$, this is
trivially true for $k\leq p$.  On the other hand, for $k=n+p$, with $n\geq1$, we have
  $$
  (\alpha g\xi)_k =
  (\alpha g\xi)_{n+p} =
  (g\xi)_n \={\ActionCocycle }
  \varphi(g,\xi|_{n-1})\xi_n  \$=
  \varphi(\g_{p+1},\xi|_{n-1})\xi_n \={$\scriptstyle \dagger$}
  \g_{n+p}\zeta_{n+q} \={ii}
  \eta_{n+p} =
  \eta_k,
  $$
  proving (a).  Focusing on (b) we have for all $n\geq1$ that
  $$
  \rho^{|\alpha|}\big(\Phi(g,\xi)\big)_{p+n} =
  \Phi(g,\xi)_n =
  \varphi(\g_{p+1},\xi|_{n-1}) \={$\scriptstyle \dagger$}
  \g_{n+p},
  $$
  proving that $\rho^{|\alpha|}\big(\Phi(g,\xi)\big) \equiv \g $, modulo $G^{(\infty)}$, hence taking care of (b).  The last two points, namely (c) and
(d) are immediate and so the proof is concluded.
  \endProof

As an immediate consequence we get a very precise description of the algebraic structure  of $\Gpd $:

\state Theorem \label ConcreteModel
  Suppose that $(\Data )$ satisfies the conditions of \lcite {\StandingHyp } and is moreover {\essfree }.  Then $\Gpd $ is
isomorphic to the sub-groupoid of \/
  $E^\infty\times (\corona \ifundef {rtimes} \times \else \rtimes \fi_{\q \rho} {\bf Z}) \times E^\infty$ given by
  $$
  \TEL = \left \{
  \matrix {
  (\eta;\q \g ,p-q;\zeta) \ \in & E^\infty\times (\corona \ifundef {rtimes} \times \else \rtimes \fi_{\q \rho} {\bf Z}) \times E^\infty : \hfill \cr
  & \g \in G^\infty,\  p,q \in {\bf N},                \hfill \pilar {12pt} \cr
  & \g_{n+p+1} = \varphi(\g_{n+p},\zeta_{n+q}),             \hfill \pilar {12pt} \cr
  & \eta_{n+p} = \g_{n+p}\zeta_{n+q}, \hbox { for all }n\geq1 \pilar {15pt}
  }
  \right \}.
  $$

\medskip
Recall from \cite {\KPRR} that the C*-algebra of every  graph is a groupoid C*-algebra for a certain groupoid constructed from
the graph, and informally called the groupoid for the \"{tail equivalence with lag}.

Viewed through the above perspective, our groupoid may also deserve such a denomination, except that the lag is not just
an integer as in \cite {\KPRR}, but an element of the lag group $\corona \ifundef {rtimes} \times \else \rtimes \fi_\rho {\bf Z}$ precisely described by the lag function
$\lag $ introduced in \lcite {\DefineLag }.

\section The topology of $\Gpd $

It is now time we look at the topological aspects of $\Gpd $.  In fact what we will do is simply  transfer the topology
of $\Gpd $ over to $\TEL $ via $F$.  Not surprisingly $F$ will turn out to be an isomorphism  of topological groupoids.

Recall from \cite [Proposition 4.14]{\actions } that, if $S$ is an inverse semigroup acting on a locally compact Hausdorff
topological space $X$, then the corresponding groupoid of germs, say $\G $, is topologized by means of the basis
consisting of sets of the form
  $$
  \Theta(s,U),
  $$
  where $s \in S$, and $U$ is an open subset of $X$, contained in the domain of the partial homeomorphism attached to $s$ by the given
action.  Each $\Theta(s,U)$ is in turn defined by
  $$
  \Theta(s,U) = \Big \{[s,x]\in\G : x\in U\Big \}.
  $$
  See \cite [4.12]{\actions } for more details.

If we restrict the choice of the $U$'s above to a predefined basis of open sets of $X$, e.g.~the collection of all
cylinders in $E^\infty$ in the present case, we evidently get the same topology on the groupoid of germs.  Therefore,
referring to the model of $\Gpd $ presented in \lcite {\FirstModel }, we see that a basis for its topology consists of the
sets of the form
  $$
  \Theta(\alpha,g,\beta;\gamma) :=   \Big \{\germ \alpha g\beta\xi \in \Gpd : \xi \in \cyl \gamma\Big \},
  \equationmark Slice
  $$
  where $(\alpha,g,\beta) \in \SGE $, and $\gamma \in E^*$.
  We may clearly suppose that $|\gamma| \geq |\beta|$ and, since   $\Theta(\alpha,g,\beta;\gamma) = \ifundef {varnothing} \emptyset \else \varnothing \fi$, unless $\beta$ is a prefix of $\gamma$, we may also
assume that $\gamma=\beta\varepsilon$, for some $\varepsilon \in E^*$.

In this case, given any $\germ \alpha g\beta\xi \in   \Theta(\alpha,g,\beta;\gamma)$, notice that $\xi\in Z(\gamma)$, and
  $$
  (\alpha,g,\beta) (\gamma,1,\gamma) =   \big(\alpha g\varepsilon,\varphi(g,\varepsilon),\gamma\big),
  $$
  from where one concludes that
  $$
  \germ \alpha g\beta\xi  = \germ {\alpha g\varepsilon}{\varphi(g,\varepsilon)}\gamma\xi,
  $$
  for all $\xi\in Z(\gamma)$, and hence also that
  $$
  \Theta(\alpha,g,\beta;\gamma) = \Theta\big(\alpha g\varepsilon,\varphi(g,\varepsilon),\gamma; \gamma\big).
  $$

  This shows that any set of the form \lcite {\Slice } coincides with another such set for which $\beta=\gamma$.  We may therefore
do away with this repetition  and redefine
  $$
  \Theta(\alpha,g,\beta) :=   \Big \{\germ \alpha g\beta\xi \in \Gpd : \xi \in \cyl \beta\Big \}.
  \equationmark ReSlice
  $$

We have therefore shown:

\state Proposition
  The collection of all sets of the form $\Theta(\alpha,g,\beta)$, where $(\alpha,g,\beta)$ range in $\SGE $, is a basis for the topology of
$\Gpd $.

We may now give a precise description of the topology of $\Gpd $, once it is viewed from the alternative point of view of
\lcite {\ConcreteModel }:

\state Proposition For each $(\alpha,g,\beta)$ in $\SGE $, the image of\/ $\Theta(\alpha,g,\beta)$ under $F$ coincides with the set
  $$
  \Omega(\alpha,g,\beta) :=
  \left \{
  \matrix {
  (\eta;\q \g ,k;\zeta) \in \TEL :
  &
    \eta \in \cyl \alpha,\
    \g \in G^\infty,\
    k = |\alpha|-|\beta|,\
    \zeta \in \cyl \beta,                                                 \hfill \cr
  & \g_{1+|\alpha|}=g,                                               \hfill \pilar {12pt} \cr
  & \g_{n+|\alpha|+1} = \varphi(\g_{n+|\alpha|},\zeta_{n+|\beta|}),                     \hfill \pilar {11pt} \cr
  & \eta_{n+|\alpha|} = \g_{n+|\alpha|}\zeta_{n+|\beta|},  \hbox { for all } n\geq1       \hfill \pilar {13pt}
  }
  \right \},
  $$
  and hence the  collection of all such sets
  form the basis for a topology on $\TEL $, with respect to which the latter  is isomorphic to $\Gpd $ as  topological groupoids.

  \Proof Left for the reader. \endProof

We may now summarize the main results obtained so far:

\state Theorem \label GroupoidModel
  Suppose that $(\Data )$ satisfies the conditions of \lcite {\StandingHyp } and is moreover {\essfree }.  Then $\OGE $ is
*-isomorphic to the C*-algebra of the groupoid $\TEL $ described in \lcite {\ConcreteModel }, once the latter is equipped with the
topology generated by the basis of open sets $\Omega(\alpha,g,\beta)$ described above, for all $(\alpha,g,\beta)$ in $\SGE $.


  \def \proj {q}	
  \def \gp {v}	    
  \def \rep {V}	    
  \def \t {t}		
  \def \modmap {\Psi}	
  \def \repalg {\psi}	
  \def \projMod {Q}	

\section $\OGE $ as a Cuntz-Pimsner algebra

Inspired by Nekrashevych's paper \cite{\NekraJO}, we will now give a description of $\OGE $ as a Cuntz-Pimsner algebra \cite{\Pimnsner}.
With this we will also be able to prove that $\OGE $ is nuclear and that $\Gpd $ is amenable when $G$ is an amenable group.
As before, we will work under the conditions of \lcite {\StandingHyp }.

We begin by introducing the algebra of coefficients over which the relevant Hilbert bimodule, also known as a
correspondence, will later be constructed.

Since the action of $G$ on $E$ preserves length \lcite {\extendedaction .\KeepLength }, we see that the set of vertexes of
$E$ is $G$-invariant, so we get an action of $G$ on $E^0$ by restriction.  By dualization $G$ acts on the algebra
$C(E^0)$ of complex valued
  functions\fn{Notice that, since $E^0$ is a finite set, $C(E^0)$ is nothing but ${\bf C}^{|E^0|}$.}
  on $E^0$.
  We may therefore form the crossed-product C*-algebra
  $$
  A = C(E^0) \ifundef {rtimes} \times \else \rtimes \fi G.
  $$

  Since $C(E^0)$ is a unital algebra, there is a canonical unitary representation of $G$ in the crossed product, which we
will denote by $\{\gp_g\}_{g\in G}$.

  On the other hand, $C(E^0)$ is also canonically isomorphic to a subalgebra of $A$ and we will therefore identify these
two algebras without further warnings.

  For each $\vr $ in $E^0$, we will denote the characteristic function of the
singleton $\{\vr \}$ by $\proj_\vr $, so that $\{\proj_\vr : \vr \in E^0\}$ is the canonical basis of $C(E^0)$, and thus $A$ coincides
with the closed linear span of the set
  $$
  \big \{\proj_\vr \gp_g: \vr \in E^0,\ g\in G\big \}.
  \equationmark LinSpanForCP
  $$

For later reference, notice that the covariance condition in the crossed product reads
  $$
  \gp_g  \proj_\vr = \proj_{g\vr }  \gp_g
  \for \vr \in E^0 \for g\in G.
  \equationmark CovarCondCrossProdinho
  $$

Our next step is to construct a correspondence over $A$.   In preparation for  this we  denote by $A^\ed $
the right ideal of $A$ generated by $\proj_{\src (\ed )} $,  for each $\ed \in E^1$.  In technical terms
  $$
  A^\ed = \proj_{\src (\ed )} A.
  $$

  With the obvious right $A$-module structure, and the inner product defined by
  $$
  \langle y,z\rangle = y^*z \for y,z\in A^\ed ,
  $$
  one has that $A^\ed $ is a right Hilbert $A$-module.  Notice that this in not necessarily a full Hilbert module since
$\langle A^\ed ,A^\ed \rangle$ is the two-sided ideal of $A$ generated by $\proj_{\src (\ed )} $, which might be a proper ideal in
some cases.

As already seen in \lcite{\LinSpanForCP}, $A$ is spanned by the elements of the form $\proj_\vr \gp_g$.  Therefore
$A^\ed$ is spanned by the elements of the form
  $
  \proj_{\src (\ed )} \proj_\vr \gp_g,
  $
  but, since the $\proj$'s are mutually orthogonal, this is either zero or equal to   $\proj_{\src (\ed )}  \gp_g$.
Therefore we see that
  $$
  A^\ed = \overline{\hbox{span}} \{\proj_{\src (\ed )} \gp_g : g\in G\}.
  $$

Introducing the right Hilbert $A$-module which will later be given the structure of a correspondence over $A$, we
define
  $$
  M = \bigoplus_{\ed \in E^1} A^e.
  $$

Observe that if $\vr $ is a vertex which is the source of many edges, say
  $$
  \src ^{-1}(\vr ) = \{\ed _1,\ed _2,\ldots,\ed_n\},
  $$
  then
  $$
  A^{\ed_i} = \proj_{\src (\ed_i)} A = \proj_\vr A,
  $$
  for all $i$, so that $\proj_\vr A$ appears many times as a direct summand of $M$.  However these copies of  $\proj_\vr
A$ should be suitably distinguished, according to which edge $\ed_i$ is being considered.

On the other hand, notice that if $\src ^{-1}(\vr )= \ifundef {varnothing} \emptyset \else \varnothing \fi$, then $\proj_\vr A$ does not appear among the summands of
$M$,  at all.

Addressing the fullness of $M$, observe that
  $$
  \langle M,M\rangle = \sum_{{\vr \in E^0 \atop \src ^{-1}(\vr )\neq\ifundef {varnothing} \emptyset \else \varnothing \fi}}A \proj_\vr A,
  $$
  so, when $E$ has no \"{sinks},
that is, when $\src ^{-1}(\vr )$ is nonempty for every $\vr $, one has that $M$ is full.

Given $\ed \in E^1$, the  element $\proj_{\src (\ed )}$, when viewed as an element of $A^\ed \subseteq M$, will play a very special
role in what follows, so we will give it a special notation, namely
  $$
  \t_\ed := \proj_{\src (\ed )}.
  \equationmark DefineTe
  $$

There is a small risk of confusion here in the sense that, if $\ed _1 ,\ed _2 \in E^1$ are such that
  $$
  \vr := \src (\ed _1)=\src (\ed _2),
  $$
  then \lcite {\DefineTe } assigns $\proj_\vr $ to both $\t_{\ed _1}$ and $\t_{\ed _2}$.  However the coordinate in which $\proj
_\vr $ appears in $\t_{\ed_i}$ is determined by the corresponding $\ed_i$, so if $\ed _1\neq\ed _2$, then $\t_{\ed _1}\neq\t_{\ed _2}$.

In order to completely dispel any confusion, here is the technical definition:
  $$
  \t_\ed = (m_\oed )_{\oed \in E^1},
  $$
  where
  $$
  m_\oed = \left \{\matrix { \proj_{\src (\ed )}, &\hbox {if } \oed =\ed ,\hfill \cr \pilar {12pt} 0, &\hbox {otherwise.}}\right .
  $$

  We should notice that
  $$
  \t_\ed \proj_{\src (\ed )} =   \t_\ed,
  \equationmark TedProj
  $$
  and  that any element of  $M$ may be written  uniquely as
  $$
  \sum_{\ed \in E^1} \t_\ed y_\ed ,
  \equationmark UniqueDecInM
  $$
  where each $y_\ed \in A^\ed $.

As the next step in  constructing a correspondence over $A$,
  we would now like to define a certain *-homomorphism from $A$ to the algebra $\Lin (M)$ of adjointable linear operators
on $M$.
  Since $A$ is a crossed product algebra, this will be accomplished once we produce a covariant representation
$(\repalg ,\rep )$ of the C*-dynamical system $\big(C(E^0),G\big)$.  We begin with the group representation $\rep$.

\definition For each $g\in G$, let $\rep_g$ be the linear operator on $M$ given by
  $$
  \rep_g\Big (\sum_{\ed \in E^1} \t_\ed y_\ed \Big ) = \sum_{\ed \in E^1} \t_{g\ed } \gp_{\varphi(g,\ed )}y_\ed ,
  $$
  whenever $y_\ed \in A^\ed $, for each $\ed $ in $E^1$.

By the uniqueness in \lcite {\UniqueDecInM }, it is clear that $\rep_g$ is well defined.

\state Proposition Each $\rep_g$ is a unitary operator in $\Lin (M)$.
Moreover, the correspondence $g\mapsto\rep_g$ is a unitary representation of $G$.

\Proof Let $g\in G$.  We begin by claiming that the expression defining $\rep_g$ above  holds true whenever the   $y_\ed $ are in $A$,
and not necessarily restricted to  $A^\ed $.  Since $\rep_g$ is clearly additive, we only need to check that
  $$
  \rep_g(\t_\ed y) = \t_{g\ed } \gp_{\varphi(g,\ed )}y
  \for y\in A.
  $$
  Observing that $\t_\ed =  \t_\ed \proj_{\src (\ed )}$,   we have
  $$
  \rep_g(\t_\ed y) =
  \rep_g(\t_\ed \proj_{\src (\ed )}y) =
  \t_{g\ed } \gp_{\varphi(g,\ed )}\proj_{\src (\ed )}y \$=
  \t_{g\ed } \proj_{\src (\varphi(g,\ed )\ed )}\gp_{\varphi(g,\ed )}y  \={\Equacoes .\PhiStarOnVert }
  \t_{g\ed } \proj_{\src (g\ed )}\gp_{\varphi(g,\ed )}y =
  \t_{g\ed } \gp_{\varphi(g,\ed )}y,
  $$
  proving the claim.  One therefore concludes that  $\rep_g$ is  right-$A$-linear.

We next claim that, for all $\ed ,\oed \in E^1$, one has
  $$
  \langle\rep_g(\t_\ed ),\t_\oed \rangle =
  \langle\t_\ed ,\rep_{g^{-1}}(\t_\oed )\rangle.
  \equationmark ClaimForAdjoint
  $$
  We have
  $$
  \langle\rep_g(\t_\ed ),\t_{\oed }\rangle =
  \langle\t_{g\ed } \gp_{\varphi(g,\ed )},\t_{\oed }\rangle =
   \gp_{\varphi(g,\ed )}^*  \langle\t_{g\ed },\t_{\oed }\rangle =
  \equal {g\ed }{\oed } \gp_{\varphi(g,\ed )}^{-1} \proj_{\src (g\ed )} \$=
  \equal {g\ed }{\oed } \proj_{\src (\varphi(g,\ed )^{-1}g\ed )} \gp_{\varphi(g,\ed )}^{-1}  \={\Equacoes .\PhiStarOnVert }
  \equal {g\ed }{\oed } \proj_{\src (\ed )} \gp_{\varphi(g,\ed )}^{-1} = (\star).
  $$
  Starting from the right-hand-side of \lcite {\ClaimForAdjoint }, we have
  $$
  \langle\t_\ed ,\rep_{g^{-1}}(\t_{\oed })\rangle =
  \langle\t_\ed ,\t_{g^{-1}\oed } \gp_{\varphi(g^{-1},\oed )}\rangle =
  \equal {e}{g^{-1}\oed } \proj_{\src (\ed )} \gp_{\varphi(g^{-1},\oed )} \$=
  \equal {ge}{\oed } \proj_{\src (\ed )} \gp_{\varphi(g,g^{-1}\oed )^{-1}} =
  \equal {ge}{\oed } \proj_{\src (\ed )} \gp_{\varphi(g,\ed )}^{-1},
  $$
  which agrees with $(\star)$ above, and hence proves claim \lcite {\ClaimForAdjoint }. If $y,z\in A$, we then have that
  $$
  \langle\rep_g(\t_\ed y ),\t_\oed z\rangle =
  y^*\langle\rep_g(\t_\ed ),\t_\oed \rangle z =
  y^*\langle\t_\ed ,\rep_{g^{-1}}(\t_\oed )\rangle z =
  \langle\t_\ed y,\rep_{g^{-1}}(\t_\oed z)\rangle,
  $$
  from where one sees that   $\langle\rep_g(\xi),\eta\rangle =  \langle\xi,\rep_{g^{-1}}(\eta)\rangle$, for all $\xi,\eta\in M$, hence proving that $\rep_g$ is an
adjointable operator with
  $
  \rep_g^* =  \rep_{g^{-1}}.
  $

Let us next prove that
  $$
  \rep_g\rep_h = \rep_{gh}
  \for g,h\in G.
  $$
  By $A$-linearity   it is  enough to prove that these operators coincide
on the set formed by the $\t_\ed $'s, which is a generating set for $M$.  We thus compute
  $$
  \rep_g\big(\rep_h(\t_\ed )\big) =
  \rep_g\big(\t_{h\ed } \gp_{\varphi(h,\ed )}\big) =
  \rep_g(\t_{h\ed }) \gp_{\varphi(h,\ed )} \$=
  \t_{gh\ed } \gp_{\varphi(g,h\ed )} \gp_{\varphi(h,\ed )}=
  \t_{gh\ed } \gp_{\varphi(gh,\ed )} =
  \rep_{gh}(\t_\ed ).
  $$

Since it is evident that $\rep_1$ is the identity operator on $M$ we obtain,
as a consequence, that
  $
  \rep_g^{-1} =
  \rep_{g^{-1}} =
  \rep_g^*,
  $
  so each $\rep_g$ is unitary and the proof is concluded.
\endProof

In order to complete our covariant pair we must now construct a *-homomorphism from $C(E^0)$ to $\Lin (M)$.  With this in
mind we give the following:

\definition For every $\vr $ in $E^0$, let
  $$
  M_\vr = \bigoplus_{\ed \in\ran ^{-1}(\vr )} A^\ed ,
  $$
  which we view as a complemented sub-module of $M$.  In addition, we let $\projMod_\vr $ be the orthogonal projection from $M$
to $M_\vr $, so that
  $$
  \projMod_\vr (\t_\ed y) = \equal {\ran (\ed )}{\vr }\t_\ed y
  \for \ed \in E^1 \for y\in A.
  \equationmark FormulaForQ
  $$

Observe that the $\projMod_\vr $ are pairwise orthogonal projections and that
  $
  \sum_{\vr \in E^0} \projMod_\vr = 1.
  $

\definition
  Let $\repalg :C(E^0) \to \Lin (M)$ be the unique unital *-homomorphism such that
  $$
  \repalg (\proj_\vr ) = \projMod_\vr
  \for \vr \in E^0.
  $$

From  our working hypothesis that $E$ has no sources, we see that  for every $\vr $ in $E^0$, there is some $\ed \in E^1$ such that
$\ran (\ed )=\vr $.  So
  $$
  \projMod_\vr (\t_\ed )=\t_\ed ,
  $$
  whence   $\projMod_\vr \neq0$.  Consequently $\repalg $ is injective.

\state Proposition The pair $(\repalg ,\rep )$ is a covariant representation of the C*-dynamical system $\big(C(E^0),G\big)$ in
$\Lin (M)$.

\Proof
All we must do is check the covariance condition
  $$
  \rep_g \repalg (y) = \repalg \big(\auto_g(y)\big) \rep_g
  \for g\in G \for y\in C(E^0),
  $$
  where $\auto $ is the name we temporarily give to the action of $G$ on $C(E^0)$.  Since $C(E^0)$ is spanned by the
$\proj_\vr $, it suffices to consider $y=\proj_\vr $, in which case the above identity becomes
  $$
  \rep_g \projMod_\vr =
  \projMod_{g\vr }  \rep_g.
  \equationmark ShortCovar
  $$

  Furthermore $M$ is generated, as an $A$-module, by the $\t_\ed $, for $\ed \in E^1$, so we only need to verify this on the
$\t_\ed $.  We have
  $$
  \rep_g \big(\projMod_\vr (\t_\ed )\big) =
  \equal {\ran (\ed )}{\vr } \rep_g (\t_\ed ) =
  \equal {\ran (\ed )}{\vr } \t_{g\ed } \gp_{\varphi(g,e)},
  $$
  while
  $$
  \projMod_{g\vr }  \big(\rep_g (\t_\ed )\big) =
  \projMod_{g\vr } \big(\t_{g\ed } \gp_{\varphi(g,e)}\big) =
  \equal {\ran (g\ed )}{g\vr }    \t_{g\ed } \gp_{\varphi(g,e)},
  $$
  verifying  \lcite {\ShortCovar } and   concluding the proof.
  \endProof

It follows from \cite [Proposition 7.6.4 and Theorem 7.6.6]{\Pedersen} that there exists a *-ho\-mo\-mor\-phism
  $$
  \modmap : C(E^0) \ifundef {rtimes} \times \else \rtimes \fi G \to \Lin (M),
  $$
  such that
  $$
  \modmap (\proj_\vr )=\projMod_\vr \for  \vr \in E^0,
  $$ and $$
  \modmap (\gp_g)=\rep_g \for g\in G.
  $$

Equipped with the left-$A$-module structure provided by $\modmap$, we then have that $M$ is a correspondence over $A$.

For later reference we record here a few useful calculations involving the left-module structure of $M$.

\state Proposition \label LeftModCalc
  Let  $g\in G$,  $\ed\in E^1$, and $\vr \in E^0$.  Then
  \Zitem $\gp_g \t_{\ed} = \t_{g\ed} \gp_{\varphi(g,\ed)}$,
  \zitem $\proj_\vr \gp_g \t_{\ed} = \equal{\ran(g\ed)}{\vr} \t_{g\ed} \gp_{\varphi(g,\ed)}.$

\Proof We have
  $$
  \gp_g  \t_{\ed}  =
  \modmap(\gp_g) \t_{\ed} =
  \rep_g (\t_{\ed} ) =
  \t_{g\ed} \gp_{\varphi(g,\ed)},
  $$
  proving (a).  Also
  $$
  \proj_\vr \gp_g \t_{\ed} =
  \modmap(\proj_\vr) (\gp_g \t_{\ed}) =
  \projMod_\vr (\t_{g\ed} \gp_{\varphi(g,\ed)})=
  \equal{\ran(g\ed)}{\vr} \t_{g\ed} \gp_{\varphi(g,\ed)}.
  \endProof

It is our next goal to prove that $\OGE$ is naturally isomorphic to the Cuntz-Pimsner algebra associated to the correspondence $M$, which we
denote by $\O_M$.  As a first step, we identify a certain Cuntz-Krieger $E$-family.

\state Proposition \label RelInOM
  The following relations hold within $\O_M$.
  \iaitem
  \aitem For every $\vr\in E^0$, one has that   $\sum_{\ed \in \ran ^{-1}(\vr )}\t_\ed \t_\ed ^* = \proj_\vr$.
  \aitem $\sum_{\ed \in E^1}\t_\ed \t_\ed ^* = 1$.
  \aitem  The set
  $
  \{\proj_\vr :  \vr \in E^0\}\cup\{\t_\ed :  \ed \in E^1\}
  $
  is a Cuntz-Krieger $E$-family.

\Proof
We first  claim that, for every  $\vr\in E^0$, and every $m\in M$, one has that
  $$
  \sum_{\ed \in \ran ^{-1}(\vr )}\t_\ed \t_\ed ^* m =   \proj_\vr m.
  $$
  To prove it, it is enough to consider the case in which $m = \t_\oed $, for $\oed\in E^1$, since these generate $M$.  In
this case we have
  $$
  \sum_{\ed \in \ran ^{-1}(\vr )}\t_\ed \t_\ed ^* \t_\oed  =
  \equal{\ran(\oed)}{\vr} \t_\oed \t_\oed ^* \t_\oed =
  \equal{\ran(\oed)}{\vr} \t_\oed \={\FormulaForQ}
  \projMod_\vr (\t_\oed) =
  \proj_\vr \t_\oed,
  $$
  proving the claim.  This says that the pair
  $
  \big(\proj_\vr,  \sum_{\ed \in \ran ^{-1}(\vr )}\t_\ed \t_\ed ^*\big)
  $
  is a redundancy or, adopting the terminology of \cite{\Pimnsner}, that the generalized compact operator
  $$
  \sum_{\ed \in \ran ^{-1}(\vr )}\Omega_{\t_\ed, \t_\ed}
  $$
  is mapped to $\modmap(\proj_\vr)$ via $\Psi^{(1)}$.
  Therefore
  $$
  \proj_\vr = \sum_{\ed \in \ran ^{-1}(\vr )}\t_\ed \t_\ed ^*,
  $$
  in $\O_M$, proving (a).  Point (b) then follows from the fact that
  $
  \sum_{\vr\in E^0} \proj_\vr =1.
  $

Focusing now on (c), it is evident that $\{\proj_\vr : \vr \in E^0\}$ is a family of mutually orthogonal projections.
Moreover, for each $e\in E^1$, we have
  $$
  \t_\ed^* \t_\ed =
  \langle\t_\ed, \t_\ed\rangle =
  \proj_{\src(\ed)},
  $$
  proving \lcite{\DefineCKFamily.\CkTwo} and also  that $\t_\ed$ is a partial isometry.  Property
\lcite{\DefineCKFamily.\CkOne} also holds in view of (a), so the proof is concluded.
\endProof

\state Proposition  \label OneMapForCuntzPimsner
  There exists a unique surjective *-homomorphism
  $$
  \Lambda:\OGE\to\O_M
  $$
  such that
  $\Lambda(p_\vr) = \proj_\vr$,
  $\Lambda(\s_\ed) = \t_\ed$, and
  $\Lambda(u_g) = \gp_g$.

\Proof
By the universal property of $\OGE$, in order to prove the existence of $\Lambda$
it is enough to check that the
  $\proj_\vr$,
  $\t_\ed$, and
  $\gp_g$
  satisfy the conditions  of \lcite{\DefineOGE}.

  Condition \lcite{\DefineOGE.a} has already been proved above while
  \lcite{\DefineOGE.b} is evidently true since $\gp$ is a representation of $G$ in $C(E^0) \ifundef {rtimes} \times \else \rtimes \fi G \subseteq \O_M$.
  Condition
  \lcite{\DefineOGE.c} is precisely           \lcite{\LeftModCalc.i}, while
  \lcite{\DefineOGE.d} was taken care of in   \lcite{\CovarCondCrossProdinho}.

Since $A$ is spanned  by the $\proj_\vr$ and the  $\gp_g$ by \lcite{\LinSpanForCP}, and since $M$ is generated over
$A$ by the $\t_\ed$, we see that $\O_M$ is spanned by the set
  $$
  \{\proj_\vr, \t_\ed, \gp_g : \vr\in E^0, \ \ed\in E^1, \ g\in G\},
  $$
  so $\Lambda$ is surjective.
  \endProof

Let us now prove that  $\Lambda$ is invertible by providing an inverse to it.  Since $A$ is the crossed product C*-algebra $C(E^0) \ifundef {rtimes} \times \else \rtimes \fi G$,
one sees that \lcite{\DefineOGE.a\&d} guarantees the existence of  a *-homomorphism
  $$
  \theta_A: A \to \OGE,
  $$
  sending the $\proj_\vr$ to the $p_\vr $, and the $\gp_g$ to the $u_g$.  For each $\ed$ in $E^1$, consider the linear mapping
  $$
  \theta_M :  M \to \OGE,
  $$
  given, for every $m = (m_\ed )_{\ed \in E^1} \in M$, by
  $$
  \theta_M(m)= \sum_{e\in E^1}\s_\ed \theta_A(m_\ed) \in \OGE.
  $$
  Notice that    $\theta_M (\t_\ed )=  \s_\ed$, for all $\ed\in E^1$, because
  $$
  \theta_M (\t_\ed )=
  \s_\ed \theta_A(\proj_{\src (\ed )}) =
  \s_\ed p_{\src (\ed )} =
  \s_\ed.
  $$

\state Lemma  The pair $(\theta_A,\theta_M)$ is a representation of the correspondence $M$ in the sense of \cite[Theorem 3.4]{\Pimnsner}, meaning that
for all $y\in A$ and all $\xi,\xi'\in M$,
  \izitem
  \zitem $\theta_M(\xi)\theta_A(y) = \theta_M(\xi y),$
  \zitem $\theta_A(y)\theta_M(\xi) = \theta_M(y\xi),$
  \zitem $\theta_M(\xi)^*\theta_M(\xi') = \theta_A(\langle\xi,\xi'\rangle).$

\Proof
  Considering the various spanning sets at our disposal, we may assume that $y\={\LinSpanForCP}\proj_\vr \gp_g$,
that $\xi = \t_{\ed} z$, and  $\xi' = \t_{\ed'} z'$,
  with $\vr \in E^0$,  $g\in G$,  $\ed,\ed'\in E^1$, $z\in\proj_{\src(\ed)} A$, and $z'\in\proj_{\src(\ed')} A$.
Then
  $$
  \theta_M(\xi)\theta_A(y) =
  \theta_M(\t_{\ed} z)\theta_A(y) =
  \s_{\ed} \theta_A(z)\theta_A(y) =
  \s_{\ed} \theta_A(zy) =
  \theta_M(\t_{\ed} zy) =
  \theta_M(\xi y),
  $$
  proving (i).  As for (ii), we have
  $$
  \theta_A(y)\theta_M(\xi) =
  \theta_A(\proj_\vr \gp_g)\theta_M(\t_{\ed} z) =
  p_\vr u_g \s_{\ed}  \theta_A(z) =
  p_\vr \s_{g\ed} u_{\varphi(g,\ed)}  \theta_A(z) \$=
  \equal{\ran(g\ed)}{\vr} \s_{g\ed} \theta_A(\gp_{\varphi(g,\ed)} z) =
  \equal{\ran(g\ed)}{\vr} \theta_M(\t_{g\ed} \gp_{\varphi(g,\ed)}z\big) \={\LeftModCalc.ii}
  \theta_M(\proj_\vr \gp_g \t_{\ed}z) =
  \theta_M(y\xi),
  $$
  proving (ii).  Focusing now on (iii), we have
  $$
  \theta_M(\xi)^*\theta_M(\xi') =   (\s_{\ed} \theta_A(z)\big)^*\s_{\ed'} \theta_A(z') =
  \equal\ed{\ed'}  \theta_A(z)^*p_{\src(e)} \theta_A(z') \$=
  \equal\ed{\ed'}  \theta_A(z^*\proj_{\src(e)} z') =
  \theta_A(\langle\xi,\xi'\rangle).
  \endProof

It is well known \cite[Theorem 3.4]{\Pimnsner} that the Toeplitz algebra for the correspondence $M$, usually denoted
$\T_M$, is universal for representations of $M$, so there exists a *-ho\-mo\-mor\-phism
  $$
  \Theta_0:\T_M\to \OGE,
  $$
  coinciding with $\theta_A$ on $A$ and with $\theta_M$ on $M$.

\state Theorem \label CPPicture
  The map $\Theta_0$, defined above, factors through $\O_M$, providing a *-isomor\-phism
  $$
  \Theta:\O_M\to \OGE,
  $$
  such that $\Theta(\proj_\vr ) = p_\vr $, $\Theta(\t_\ed )=\s_\ed $, and $\Theta(\gp_g) = u_g$, for all $\vr\in E^0$, $\ed\in E^1$, and $g\in G$.

\Proof
  The factorization property
  follows immediately from \lcite{\RelInOM.b} and an easy modification of \cite[Proposition 7.1]{\ExelVesshik} to
Cuntz-Pimsner algebras.

In order to prove that $\Theta$ is an isomorphism, observe that $\Theta\circ\Lambda$
coincides with the identity map on the generators of $\OGE$, by \lcite{\OneMapForCuntzPimsner}, and hence $\Theta\circ\Lambda=id$.  The
result then follows from the fact that $\Lambda$ is surjective.
  \endProof

\state Corollary \label Amenabuilidade
  If $G$ is amenable then $\OGE$ is nuclear.

\Proof
  The amenability of $G$ ensures that $C(E^0) \ifundef {rtimes} \times \else \rtimes \fi G$ is nuclear.  The result then follows from
\lcite{\CPPicture},  the fact that Toeplitz-Pimsner algebras over nuclear coefficient algebras is nuclear
\cite[Theorem 4.6.25]{\BO}, and so are quotients of nuclear algebras \cite[Theorem 9.4.4]{\BO}.
\endProof

\state Remark \rm Since $E^0$ is finite, the nuclearity of  $C(E^0) \ifundef {rtimes} \times \else \rtimes \fi G$ is equivalent to the amenability of $G$.
However, if the present construction is generalized for infinite graphs, one could produce examples of non amenable
groups acting amenably on $E^0$, in which case $C(E^0) \ifundef {rtimes} \times \else \rtimes \fi G$ would be amenable.  The proof of  \lcite{\Amenabuilidade}
could then
be adapted to prove that $\OGE$ is nuclear.

\state Corollary If $G$ is amenable and $(\Data)$ is {\essfree} then $\Gpd$ and its sibling  $\TEL$ are  amenable
groupoids.

\Proof Follows from \lcite{\Amenabuilidade}, \lcite{\GroupoidModel}, and \cite[Theorem 5.6.18]{\BO}.
 \endProof

Nekrashevych has proven in \cite[Theorem 5.6]{\NC}, that a certain groupoid of germs, denoted ${\cal D}_G$, constructed
in the context of self-similar groups, is amenable under the hypothesis that the group is \"{contracting} and
\"{self-replicating}.  Even though there are numerous differences between ${\cal D}_G$ and $\TEL$, including a different
notion of \"{germs} and Nekrashevych's requirement that group actions be \"{faithful}, we believe it should be
interesting to try to generalize Nekrashevych's result to our context.

\endgroup

\bigskip \noindent \tensc
  Departamento de Matem\'atica;
  Universidade Federal de Santa Catarina;
  88010-970 Florian\'opolis SC;
  Brazil
  \hfill\break \tt (ruyexel@gmail.com)

\bigskip \noindent \tensc
  {Departamento de Matem\'aticas, Facultad de Ciencias;
  Universidad de C\'adiz, Campus de Puerto Real;
  11510 Puerto Real (C\'adiz);
  Spain
  \hfill\break \tt (enrique.pardo@uca.es)

\bye
\bye